\documentclass[twoside,12pt]{article}

\usepackage{amssymb,amsmath,verbatim,bm,color,mathtools,multicol}
\usepackage{graphicx}
\usepackage{psfrag}
\usepackage{indentfirst}
\usepackage{hyperref} 
\usepackage{newtxtext}
\usepackage{algorithm}
\usepackage[noend]{algpseudocode}
\usepackage{tabularx}
\usepackage{multirow,array,booktabs}
\usepackage{caption}
\usepackage{subcaption}
\usepackage{wrapfig}
\usepackage{textcomp}
\usepackage{ulem}

\usepackage{graphicx}

\DeclareMathOperator\erf{erf}

\usepackage{tikz}
\usetikzlibrary{fit,positioning,automata,backgrounds}

\usepackage[abbrvbib,preprint,nohyperref]{jmlr2e}

\newcommand{\vectorstyle}[1]{\boldsymbol{\mathbf{#1}}}

\newcommand{\matrixstyle}[1]{\mathrm{#1}}

\firstpageno{1}
\begin{document}

\title{Differential Flatness of Slider-Pusher Systems for Constrained Time Optimal Collision Free Path Planning}
\author{Tom Lefebvre\textsuperscript{1,2}, Sander De Witte\textsuperscript{1,2}, Thomas Neve\textsuperscript{1,2} \& Guillaume Crevecoeur\textsuperscript{1,2}\\
	\textsuperscript{1}Department of Electromechanical, Systems and Metal Engineering, Ghent University, Technologiepark 131, Ghent, Belgium.\\
	\textsuperscript{2}MIRO corelab, Flanders Make, Belgium.}


\editor{}

\maketitle

\begin{abstract}
			In this work we show that the differential kinematics of slider-pusher systems {are differentially flat} assuming quasi-static behaviour and frictionless contact. Second we demonstrate that the state trajectories are invariant to time-differential transformations of the path parametrizing coordinate. For one this property allows to impose arbitrary velocity profiles on the slider without impacting the geometry of the state trajectory. This property implies that certain path planning problems may be decomposed approximately into a strictly geometric path planning and an auxiliary throughput speed optimization problem. Building on these insights we elaborate a numerical approach tailored to constrained time optimal collision free path planning and apply it to the slider-pusher system.
\end{abstract}

	\section{Introduction}
The ability to manipulate an object by pushing it with another object, is a resourceful skill for robotic systems to master. Adopting the principle of pushing can enable a system to manipulate objects that would otherwise be too large, too heavy, or too cluttered to be grasped \cite{yu2016more}. Furthermore, we argue it can significantly reduce the requirements for the end-effector's design without severely impacting the manoeuvrability of the object. 

However, manipulating an object by pushing it comes at a cost given that controlling the system becomes complicated. {The origin of the control challenge roots back to the model \cite{hogan2020feedback}}. The mechanics at the contact interface of the interacting objects become significantly more involved. The dynamics of pushing, particularly the dynamics of the combined slider-pusher system have been the subject of numerous studies. The earliest references in the mechanics and robotics community date back almost three decades \cite{mason1986mechanics,goyal1991planar1,goyal1991planar2} with recent literature focussing on data-driven probabilistic model approaches \cite{yu2016more,bauza2017probabilistic}. Though, full descriptions of the Newtonian mechanics of the system, meaning full modelling of the interaction forces, are rarely considered. A more elegant modelling approach is to adopt the quasi-static assumption which implies that the slider-pusher motions are slow enough that inertial forces are negligible compared to frictional forces. As a result a differential kinematic model is obtained relating the input velocity of the pusher to the velocities of the planar configuration of the slider.

Still the combined slider-pusher remains a hybrid system, exhibiting different contact modes depending on the actual actuation regime. It also remains underactuated system, considering that the planar configuration of the object is controlled through a single (idealized) contact point \cite{hogan2020feedback}. At the result of these mathematical system properties, manipulation of a slider through an actuated pusher requires the practice of planning and other advanced model-based predictive control approaches, which imposes significant requirements on the computational resources of the control system. Presumably that is why slider-pusher systems have only recently resurfaced in the literature 	\cite{doshi2020hybrid,raghunathan2022pyrobocop,hogan2020feedback} with some exceptions dating back longer \cite{lynch1992manipulation}.

The motivating question of this contribution is the following. What happens when we assume that the contact between slider and pusher is frictionless or at least that the friction forces between the slider and pusher are negligible compared to the friction forces between the slider and the supporting surface? These modelling conditions have been investigated for the slider-pusher system by e.g. \cite{lynch1992manipulation} but the implications with regard to control were never recognized. In particular, in this contribution we show that under this assumption the slider-pusher is differentially flat \cite{fliess1995flatness}.

Differential flatness is a structural property {that certain} nonlinear dynamical systems {exhibit }\cite{rigatos2015differential}. If a system is \textit{flat}, this denotes that all defining and differentially dependent system variables (states and controls) can be written in terms of a specific set of differentially independent variables and their derivatives \cite{rigatos2015differential}. Flatness is a resourceful property for both the analysis and controller synthesis of nonlinear dynamical systems. It is particularly advantageous for solving trajectory planning problems and trajectory tracking. Stoican et al. developed a great deal of theory to address path planning problems exploiting differential flatness and B-splines path parametrizations \cite{stoican2015flat,stoical2016obstacle,stoican2017constrained}. Flatness is also used abundantly in the design and synthesis of asymptotical set-point following control \cite{greeff2018flatness,helling2020flatness,faessler2017differential,aguilar2012trajectory}. Most of these works are tailored to Unmanned (Aerial and Ground) Vehicles (UAV and UGV), provided that the kinematic car and quadcopter dynamics are flat, though some works have extended flatness to other practical areas such as electromagnetic actuation systems \cite{thounthong2018nonlinear,jovsevski2015flatness}. In recent work Greeff et al. discuss a robust adaptive control strategy for quadcopters using flatness and the theory of Gaussian processes \cite{greeff2020exploiting}.

Unfortunately there exists no straightforward procedure to verify whether a set of dynamic system equations is differentially flat \cite{rigatos2015differential}. It remains therefore standard practice to verify for each system independently. Known systems that are flat are kinematic cars, quadcopters, gantry cranes, etc. In this work we show that the slider-pusher system can be added to this list. {Apart from showing that the slider-pusher system is differentially flat, we show that the system exhibits an auxiliary invariance property. The geometry of the state trajectory is invariant to the relation between the path coordinate and time. This property has interesting consequences. In this work we show how the property allows to decompose time optimal path planning problems in an approximate manner. }

The contributions of this work are threefold 
\begin{enumerate}
	\item First, we provide an original derivation for the differential kinematics of slider-pushers with negligible contact friction.
	\item Second, we show that the model is differentially flat and invariant to differential transformations of the path coordinate.
	\item Finally, we exploit these properties to develop a two step procedure tailored to constrained time optimal path planning which extends trivially to other path invariant flat systems.
\end{enumerate}

\section{Quasi-static model}\label{sec:quasi-static-model} In this section we propose an original derivation of the quasi-static model for slider-pusher systems. The model is physically valid in the particular case where 
\begin{enumerate}
	\item the pusher motions are slow enough that inertial forces are negligible compared to frictional forces
	\item the friction forces at the contact point are negligible with respect to the friction forces between the slider and the ground, equivalently we assume the local contact friction coefficient is equal to zero \cite{lynch1992manipulation}
\end{enumerate}

\subsection{Kinematics} 
\begin{figure}
	\centering
	\includegraphics[width=.7\columnwidth]{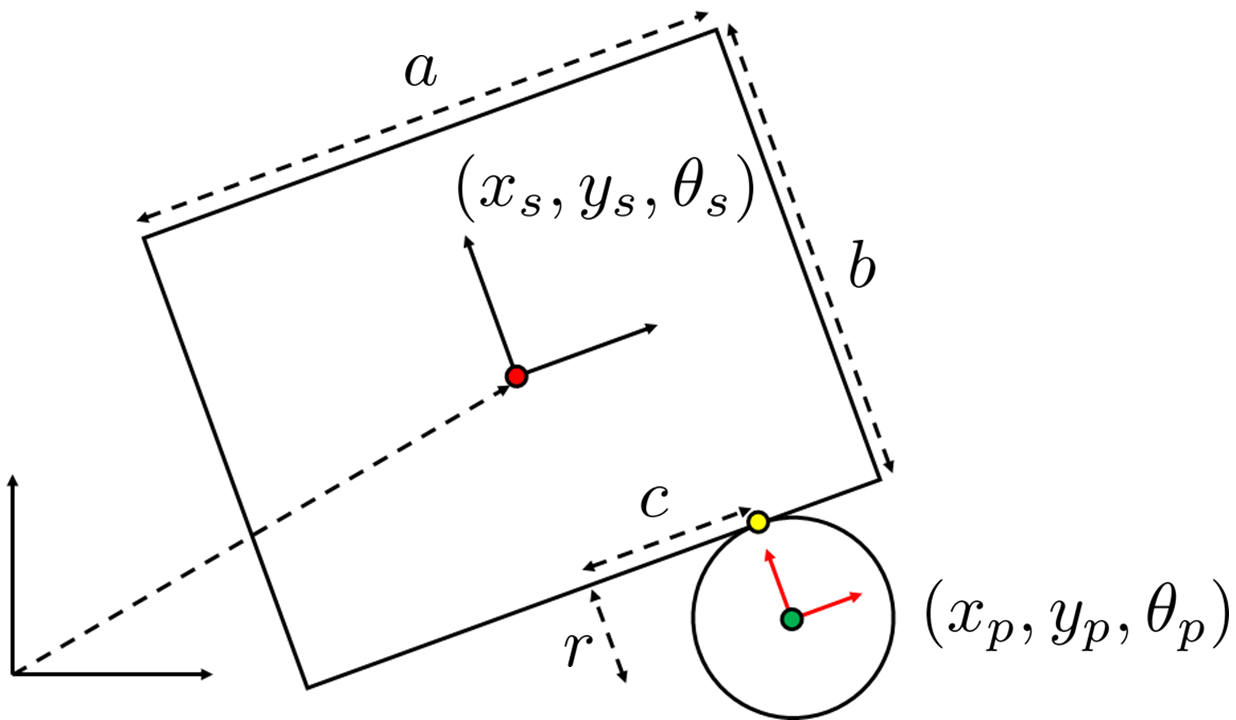}
	\caption{Kinematics of the slider-pusher.}
\label{fig:sp}
\end{figure}

Consider the schematic representation of the slider-pusher system in Fig. \ref{fig:sp}. The slider is a planar object. As such its configuration can be parametrized by its Cartesian coordinates, $x_s$, and, $y_s$, and its planar orientation, $\theta_s$, expressed in a global frame of reference. Note that in the present study, we consider rectangular sliders. The dimensions of the slider are thus quantified by its width, $a$, and its length, $b$. Second we assume that contact is maintained during manoeuvring and that the pusher has negligible dimensions (for now). Then, the position of the pusher can be defined through the position of the contact point relative to the planar pose of the slider. The position of the contact point relative to the width symmetry axis of the slider is denoted as $c$. Note that this parametrization is arbitrary so that the same model can be used to model the situation where for example the pusher switches sides.

The state, $\vectorstyle{x}_s\in \mathbb{R}^4$, of the \textit{slider-point} system is given by
\begin{equation}
	\vectorstyle{x}_s = \begin{pmatrix}
		x_s \\ y_s \\ \theta_s \\ c
	\end{pmatrix}
\end{equation}

{For irregular but smooth slider circumferences, the contact point of the slider could also be parametrized with an angle, $\phi$, and, some function $r(\phi)$. Here $\phi$ quantifies the angle between the present contact point and some arbitrary reference with respect to the centre of mass and $r(\phi)$ denotes the distance w.r.t. the centre of mass. For the sake of simplicity here we only consider rectangular sliders.}

It is assumed that the velocity of the pusher expressed in the local frame of reference, i.e. that attached to the slider, can be controlled directly. Thus, the control variable, $\vectorstyle{u}_s\in\mathbb{R}^2$, is defined as 
\begin{equation}
	\vectorstyle{u}_s = \begin{pmatrix}
		v_t \\ v_n
	\end{pmatrix}
\end{equation}
where $v_t$ and $v_n$ denote the tangential and normal velocity respectively {(indicated in red in Fig. \ref{fig:sp})}.

Finally we extend the \textit{slider-point} system to a \textit{slider-pusher} system. We model the pusher as a sphere with radius $r$. It is further assumed that the pusher is realised as a kinematic car which might be the case in practical applications. The pusher's state, $\vectorstyle{x}_p \in\mathbb{R}^3$, is then given by the planar pose of the pusher
\begin{equation}
	\vectorstyle{x}_p = \begin{pmatrix}
		x_p \\ y_p \\ \theta_p
	\end{pmatrix}
\end{equation}

The control input, $\vectorstyle{u}_p\in\mathbb{R}^2$ of the kinematic car, and hence of the pusher, are given by its local normal and angular velocity
\begin{equation}
	\vectorstyle{u}_p = \begin{pmatrix}
		v_p \\ \omega_p
	\end{pmatrix}
\end{equation}

\subsection{Quasi-static model with frictionless contact}

The quasi-static differential kinematics of the slider-{pusher} system with frictionless contact are governed by the following set of equations. {Here the parameter $\beta$ denotes a geometrical factor that will be made explicit later.}
\begin{equation}
	\label{eq:dk}
	\begin{aligned}
		\dot{x}_s &= - \frac{\beta^2}{\beta^2 + c^2}v_n \sin(\theta_s)\\
		\dot{y}_s &=  \frac{\beta^2}{\beta^2 + c^2}v_n \cos(\theta_s) \\
		{\dot{\theta}_s} &= \frac{c}{\beta^2 + c^2 }v_n \\
		\dot{c} &= v_t - \left(\frac{b}{2}+r\right) \frac{c }{\beta^2 + c^2} v_n 
	\end{aligned}
\end{equation}

We discuss two strategies to arrive at these differential kinematics. The first derives from the principle of least work. As a result we do not explicitly model the contact forces. The second approach is adopted from literature and is based on the concept of a \textit{limit surface} \cite{goyal1991planar1,lynch1992manipulation,hogan2020feedback}. Both derivation strategies pursue a quasi-static model.  The motion should further comply with the kinematic constraint imposed by the frictionless contact. 

\subsubsection{Principle of least work}\label{sec:principle-of-least-work}

Without modelling the friction forces explicitly, it is difficult to see how the slider will react when a velocity is imposed to the contact point. To remedy this issue, we will rely on the tendency of physical systems to follow paths that are associated to the least amount of work. 
\begin{enumerate}
\item Therefore we first express the energy dissipation rate, $\dot{W}$, at an infinitesimal surface element of the slider as a function of its differential state. To that end we assume that the local dissipation rate is equal to $\text{d}\dot{W} = \dot{\vectorstyle{p}}\cdot \text{d}\vectorstyle{f}$, where $\dot{\vectorstyle{p}}$ denotes the velocity, and, $\text{d}\vectorstyle{f}$ an infinitesimal local contribution to the total friction force.
\item Second we assume that the direction of $\text{d}\vectorstyle{f}$ is parallel to the local velocity $\dot{\vectorstyle{p}}$ and that the amplitude of $\text{d}\vectorstyle{f}$ is proportional to the local contribution to the total normal force, $F = m g$, { and, proportional to the local velocity according a viscous friction law}. As a result $\text{d}\vectorstyle{f}\propto \rho(\vectorstyle{q}) \text{d}A$ where $\rho(\vectorstyle{q})$ denotes the pressure distribution, $\vectorstyle{q}$ is local coordinate relative to the object's centre of mass and $\text{d}\vectorstyle{A}$ denotes an infinitesimal surface element. Further remark that by definition $F=\int \rho(\vectorstyle{q})\text{d}A$. {In the end we have that $\text{d}\dot{W} \propto \rho(\vectorstyle{q}) \|\dot{\vectorstyle{p}}\|^2$.}
\item It is not straightforward how to determine this distribution. Therefore we will assume that the contribution of each surface element is equally important. Equivalently, we assume that the pressure distribution, $\rho(\vectorstyle{q})$, is uniform. 
\item Finally the motion of the contact point should also comply with that of the pusher. The differential motion of the slider is then postulated to minimize the total dissipation of energy. 
\end{enumerate}

For frictionless contact this approach gives rise to the following constrained optimization problem. 
\begin{equation}
\label{eq:prob1}
\begin{aligned}
	\min_{\dot{\vectorstyle{x}}_s} ~& J(\dot{\vectorstyle{x}}_s;{\vectorstyle{x}}_s,{\vectorstyle{u}}_s) \\
	\text{ s.t. } &\dot{\vectorstyle{p}}_c^g = \mathrm{R}(\theta_s) \vectorstyle{u}_s
\end{aligned}
\end{equation}
Here 
\begin{enumerate}
\item $\vectorstyle{p}_c^g$ is the position of the contact point in global coordinates
$$\vectorstyle{p}_c^g = \matrixstyle{R}({\theta_s}) \vectorstyle{p}_c^l +\vectorstyle{p}_s^g$$
\item $\matrixstyle{R}(\cdot)$ is the planar rotation matrix
\item $\vectorstyle{p}_c^l$ denotes the local position of the contact point
$$\vectorstyle{p}_c^l = \begin{pmatrix}
	c \\ -\frac{b}{2}-r
\end{pmatrix}$$
\item and $\vectorstyle{p}_s^g = (x_s,y_s)$ denotes the global position of the slider.
\end{enumerate}

In agreement with the modelling assumptions listed above, the objective can be defined as
\begin{equation}
\label{eq:objective}
J(\dot{\vectorstyle{x}}_s;{\vectorstyle{x}}_s,{\vectorstyle{u}}_s) =  \int_{\mathcal{A}}  \|\dot{\vectorstyle{p}}\|^2 \rho(\vectorstyle{q})\text{d}A 
\end{equation} 
Here $\vectorstyle{p}$ represents a point on the surface in global coordinates
$$\vectorstyle{p} =\matrixstyle{R}(\theta_s) \vectorstyle{q} + \vectorstyle{p}_s^g$$
so that 
\begin{equation}
\label{eq:obj1}
J(\dot{\vectorstyle{x}}_s;{\vectorstyle{x}}_s,{\vectorstyle{u}}_s)  \propto (\dot{x}_s^2+\dot{y}_s^2) + \beta_1^2 \dot{\theta}_s^2 
\end{equation}
Solving (\ref{eq:prob1}), then yields the differential kinematic model from (\ref{eq:dk}). 

The solution still depends on the geometric factor $\beta_1^2$. This term is defined as
\begin{equation}
\beta_1^2 = \frac{1}{F}\int_{\mathcal{A}} \|\vectorstyle{q}\|^2 \rho(\vectorstyle{q})\text{d}A 
\end{equation}

In retrospect, it now is useful to analyse the objective in (\ref{eq:obj1}). There are two terms that contribute to the objective. The first term is directly related to the linear motion of the object. The second term is directly related to the angular motion of the object. It follows that the relative contribution of these two terms will affect the behavioural tendencies of the slider-pusher interaction. Now, let us first note that the linear term is unaffected by the pressure distribution, $\rho(\vectorstyle{q})$. The contribution of the angular motion on the other hand is determined by the geometric factor, $\beta_1^2$, which in turn does depend on the pressure distribution, $\rho(\vectorstyle{q})$. As a result, the pressure distribution will determine the relative contribution of the two terms and thus the behavioural tendencies of the object. If the object’s support pressure is concentrated near the boundaries of the object, it is in the object's interest to resist angular motion. Equivalently, the factor, $\beta_1^2$, will be large and so will be the angular motion's contribution to the objective. Likewise, if the object's support pressure is concentrated near the centre of mass, it is in the object's best interest to resist linear motion; the factor, $\beta_1^2$, will be small and so will be the linear motion's contribution.

For a uniform pressure distribution and a rectangular object the factor, $\beta_1^2$, equals
\begin{equation}
\beta_1^2 = \frac{1}{A}\int_{-\frac{a}{2}}^{\frac{a}{2}}\int_{-\frac{b}{2}}^{\frac{b}{2}} \|\vectorstyle{q}\|^2 \text{d}q_x\text{d}q_y = \tfrac{1}{12}D^2 
e	\end{equation}
where $D=\sqrt{a^2+b^2}$.

\begin{figure}
\centering
\includegraphics[width=\columnwidth]{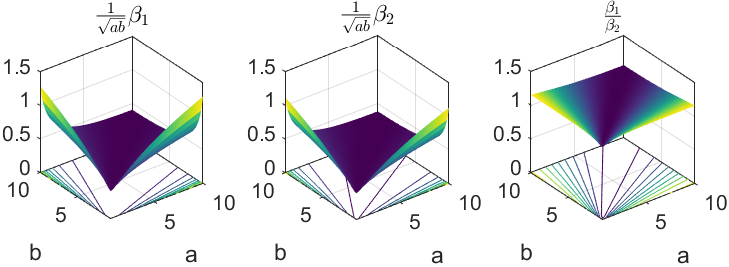}
\caption{Analytical comparison of $\beta_1$ and $\beta_2$.}
\label{fig:beta}
\end{figure}

\begin{figure*}[h]
\centering
\includegraphics[width=.495\columnwidth]{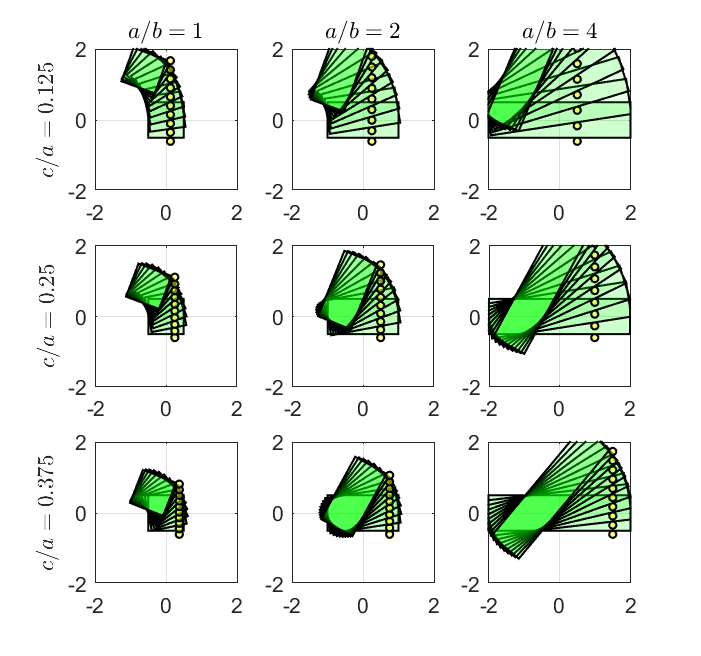}
\includegraphics[width=.495\columnwidth]{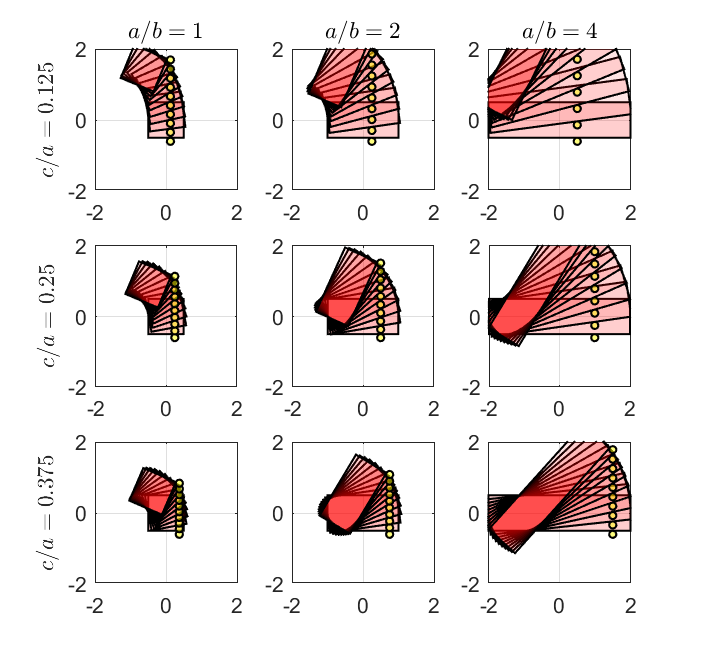}
\caption{Behavioural comparison of $\beta_1$ (left) and $\beta_2$ (right).}
\label{fig:ol}
\end{figure*}

\subsubsection{Limit surfaces}
The second derivation strategy relies on the concept of a \textit{limit surface} and is adopted from \cite{goyal1991planar1,lynch1992manipulation,hogan2020feedback}. 

The limit surface of a sliding object is defined as a closed convex surface in the space $(f_x,f_y,n)$ that encloses the origin. Here $f_x$, $f_y$ and $n$ denote the horizontal and vertical forces and planar moment that are exerted on the slider. The quasi-static assumption implies that any force exerted to the slider that does provoke a velocity response but does not provoke an acceleration, must be on the slider's limit surface. As in \cite{lynch1992manipulation}, we approximate the slider's limit surface with an ellipsoidal surface
\begin{equation}
\label{eq:limit}
\left(\frac{f_x}{f_x^*}\right)^2 + \left(\frac{f_y}{f_y^*}\right)^2 + \left(\frac{n}{n^*}\right)^2  = 1
\end{equation}
Here $f^*_x$, $f_y^*$ and $n^*$ are the maximal frictional forces and moment.
\begin{enumerate}
\item The maximal forces $f^*_x$ or $f_y^*$ are determined as the maximal linear friction force and are therefore proportional to the normal force $f_n$. It follows that $f_x^*$ and $f_y^*$ are unaffected by the pressure distribution $\rho(\vectorstyle{q})$.
\begin{equation}
	f_x^* = f_y^* = \mu \int \rho(\vectorstyle{q}) \text{d}A = \mu F
\end{equation}
\item The maximal planar moment is determined as
\begin{equation}
	n^* = \mu \int \|\vectorstyle{q}\|\rho(\vectorstyle{q}) \text{d}A = \mu F \beta_2
\end{equation}
\end{enumerate} 

Based on the limit surface, the slider's motion is further constrained until a system of equations is encountered that can be solved for $\dot{\vectorstyle{x}}_s$.
\begin{enumerate}
\item It is further assumed that the reciprocal motion, $(v_x,v_y,\omega)$ (defined at the centre of mass), of the slider is normal to the limit surface. Together with the limit surface approximation (\ref{eq:limit}), this deprives the problem of two degrees of freedom.
\item Second, the motion of the contact {and} slider are related.
\item Third, the force exerted on the slider by the pusher must go through the contact, posing an additional constraint on $n$.
\end{enumerate}
To complete the model, the velocity of the contact point is related to the velocity of the pusher. For additional details we refer to \cite{lynch1992manipulation,hogan2020feedback}. 

For frictionless contact, ergo $\mu_c=0$, the situation that is assumed here throughout, the differential kinematics collapse onto the same expressions given in (\ref{eq:dk}) but with an alternative geometric factor $\beta$.
\begin{equation}
\beta_2 = \frac{1}{F}\int_{\mathcal{A}} \|\vectorstyle{q}\| \rho(\vectorstyle{q})\text{d}A
\end{equation}
For a uniform pressure distribution and a rectangular object the factor, $\beta_2$, equals (see appendix \ref{sec:derivation-of-beta2})
\begin{equation}
\beta_2 = \frac{a^2}{12b} \log \frac{D + b}{a} -\frac{b^2}{12a}\log \frac{D - a}{b} + \frac{1}{6}D
\end{equation}

\subsubsection{Comparison}

As it turns out, the models are equivalent up to the geometric factor, $\beta$.

An analytical comparison of the geometric factors $\beta_1$ and $\beta_2$ is made in Fig. \ref{fig:beta}. The comparison is made as a function of the dimensions $a$ and $b$, the factors are scaled by their geometric mean to make the factors dimensionless. Clearly, the factors are almost equivalent. Their ratio is approximately constant but $\beta_1$ is larger than $\beta_2$ by roughly $20\%$. 

Based on the analytical comparison from Fig. \ref{fig:beta} it remains difficult to asses{s} which factor will result in the most physically accurate behaviour. The behavioural difference is however negligible as can be seen in Fig. \ref{fig:ol}. Here a behavioural comparison is {made} between both models for open-loop pushing for varying ratios for $a$ and $b$ and different initial points of contact, $c$. The pusher moves vertically and the simulation is stopped once contact is lost. The velocity of the pusher is constant. Its magnitude is irrelevant for the outcome of the numerical experiment due to a particular invariance property of the slider-pusher. For details we refer to section \ref{sec:imposing-arbitrary-velocity-profiles}.

For an analysis of the pressure distributions in function of the geometric factors, we refer to appendix \ref{sec:reciprocal-pressure-distributions}.

{\subsection{General slider geometries} The analysis from section \ref{sec:principle-of-least-work} can be repeated for an arbitrary slider geometry resulting into a generalised differential kinematic model. Such a derivation is included in Appendix \ref{sec:generalised-slider-geometry}. For simplicity, in this work we focus on sliders with a rectangular geometry. }

\section{Differential flatness}

In this section we show that the {rectangular} slider-pusher system with frictionless contact is a differentially flat system.

Let us recall the formal definition of differential flatness \cite{fliess1995flatness}. 
\begin{definition}
The system, $\dot{\vectorstyle{x}}=\vectorstyle{f}(\vectorstyle{x},\vectorstyle{u})$, with state $\vectorstyle{x}\in\mathcal{X}\subset\mathbb{R}^{n_x}$ and input $\vectorstyle{u}\in\mathcal{U}\subset\mathbb{R}^{n_u}$, is differentially flat if there exists a variable $\boldsymbol{\zeta}\in\mathcal{Z}\subset\mathbb{R}^{n_\zeta}$, whose components are differentially independent, and operators $\Lambda$, $\Phi$ and $\Psi$ such that the following holds \cite{fliess1995flatness}:

\begin{equation}
	\begin{aligned}
		\boldsymbol{\zeta} &= \Lambda(\vectorstyle{x},\vectorstyle{u},\dot{\vectorstyle{u}},\dots,\vectorstyle{u}^{(\lambda)})\\
		\vectorstyle{x} &= \Phi(\vectorstyle{\zeta},\dot{\vectorstyle{\zeta}},\ddot{\vectorstyle{\zeta}},\dots,\vectorstyle{\zeta}^{(\phi)})\\
		\vectorstyle{u} &= \Psi(\vectorstyle{\zeta},\dot{\vectorstyle{\zeta}},\ddot{\vectorstyle{\zeta}},\dots,\vectorstyle{\zeta}^{(\phi-1)})
	\end{aligned}
\end{equation}
Here $\Lambda$, $\Phi$ and $\Psi$ are smooth function operators, $\lambda$ and $\phi$ are the maximum orders of the derivatives of $\vectorstyle{u}$ and $\vectorstyle{\zeta}$ needed to describe the system and $\vectorstyle{\zeta}$ is called the flat coordinates.
\end{definition}

Intuitively, the flat coordinate $\vectorstyle{\zeta}(t)$ can be interpreted as a minimal dynamical representation of any feasible state-action trajectory, $(\vectorstyle{x},\vectorstyle{u})(t)\in\mathcal{F}\subset\mathcal{X}\times\mathcal{U}$ of the system. The feasible state-action function space, $\mathcal{F}$, is a subspace from the function space $\mathcal{X}\times\mathcal{U}$ so that any function element satisfies the dynamic constraint $\dot{\vectorstyle{x}}=\vectorstyle{f}(\vectorstyle{x},\vectorstyle{u})$.
\begin{equation}
\label{eq:feasible}
\mathcal{F} = \{(\vectorstyle{x},\vectorstyle{u})|\vectorstyle{x}\in\mathcal{X},\vectorstyle{u}\in\mathcal{U}:\dot{\vectorstyle{x}}=\vectorstyle{f}(\vectorstyle{x},\vectorstyle{u})\}
\end{equation}

Consequently, we can think of $\Lambda$ as a \textit{projection} operator from the feasible state-action function space, $\mathcal{F}$, to the motion's minimal representation in flat space, $\mathcal{Z}$. Analogously, we can think of $(\Phi,\Psi)$ as an \textit{inflation} operator from the flat function space, $\mathcal{Z}$, to the feasible state-action space, $\mathcal{F}$. It is interesting to remark that the relation is (usually) bijective. Any element in $\mathcal{F}$ is associated to an element in $\mathcal{Z}$ and vice versa \cite{fliess1995flatness}.

Systems that are known to be flat are quadcopters \cite{faessler2017differential}, gantry cranes \cite{fliess1995flatness}, cars with trailers, ... but also fully actuated multi-body systems, considering the flat coordinate $\vectorstyle{\zeta}=\vectorstyle{q}$ with generalised coordinates $\vectorstyle{q}$ and $\vectorstyle{\tau}(\vectorstyle{q},\dot{\vectorstyle{q}},\ddot{\vectorstyle{q}})=\matrixstyle{M}(\vectorstyle{q})\ddot{\vectorstyle{q}}+\vectorstyle{c}(\vectorstyle{q},\dot{\vectorstyle{q}})$. 

Here we argue that also the quasi-static slider-pusher with frictionless contact is flat.
\begin{theorem}
\label{th:1}
The quasi-static slider-pusher system with frictionless contact, (\ref{fig:sp}), is differentially flat with flat coordinates
\begin{equation}
	\vectorstyle{\zeta} = \vectorstyle{p}_s^g = \begin{pmatrix}
		x_s\\y_s
	\end{pmatrix}
\end{equation}
\end{theorem}

The flat expressions for the slider are given by 
\begin{equation}
\begin{aligned}
	&\theta_s = -\arctan \tfrac{\dot{x}_s}{\dot{y}_s} \\
	&c = \beta^2 \tfrac{\dot{x}_s\ddot{y}_s - \ddot{x}_s\dot{y}_s}{\sqrt{\dot{x}_s^2+\dot{y}_s^2}^{3}} \\
	&v_t = \alpha\tfrac{\dot{x}_s\ddot{y}_s-\ddot{x}_s\dot{y}_s}{\dot{x}_s^2+\dot{y}_s^2}+\beta^2\tfrac{\dot{x}_s\dddot{y}_s - \dddot{x}_s\dot{y}_s}{\sqrt{\dot{x}_s^2+\dot{y}_s^2}^{3}} + 3\beta^2 \tfrac{\left(\ddot{x}_s\dot{y}_s-\dot{x}_s\ddot{y}_s\right)\left(\dot{x}_s\ddot{x}_s+\dot{y}_s\ddot{y}_s\right)}{\sqrt{\dot{x}_s^2+\dot{y}_s^2}^{5}} \\
	&v_n = \left(1+\beta^2\tfrac{\left(\dot{x}_s\ddot{y}_s-\ddot{x}_s\dot{y}_s\right)^2}{\left(\dot{x}_s^2+\dot{y}_s^2\right)^3}\right)\sqrt{\dot{x}_s^2+\dot{y}_s^2} 
\end{aligned}
\end{equation}
where $\alpha=\frac{b}{2}+r$ is introduced for notational convenience.

For the pusher we retrieve the following flat expressions
\begin{equation}
\begin{aligned}
	x_p &= x_s +\alpha\sin\theta_s + c\cos\theta_s  \\
	y_p &= y_s - \alpha\cos\theta_s + c \sin\theta_s  \\
	\theta_p &= -\arctan \tfrac{\dot{x}_p}{\dot{y}_p} = f_{\theta}\left(\dot{\vectorstyle{\zeta}},\ddot{\vectorstyle{\zeta}},\dddot{\vectorstyle{\zeta}}\right)\\
	v_p &= \sqrt{\dot{x}_p^2+\dot{y}_p^2} = f_{v}\left(\dot{\vectorstyle{\zeta}},\ddot{\vectorstyle{\zeta}},\dddot{\vectorstyle{\zeta}}\right)\\
	\omega_p &= \tfrac{\dot{x}_p \ddot{y}_p - \ddot{x}_p\dot{y}_p}{\dot{x}_p^2 + \dot{y}_p^2} = f_{\omega}\left(\dot{\vectorstyle{\zeta}},\ddot{\vectorstyle{\zeta}},\dddot{\vectorstyle{\zeta}},\ddddot{\vectorstyle{\zeta}}\right)
\end{aligned}
\end{equation}

Details regarding a derivation can be found in appendix \ref{sec:derivation-of-differential-flatness}.

\begin{figure}
\centering
\includegraphics[width=.5\columnwidth]{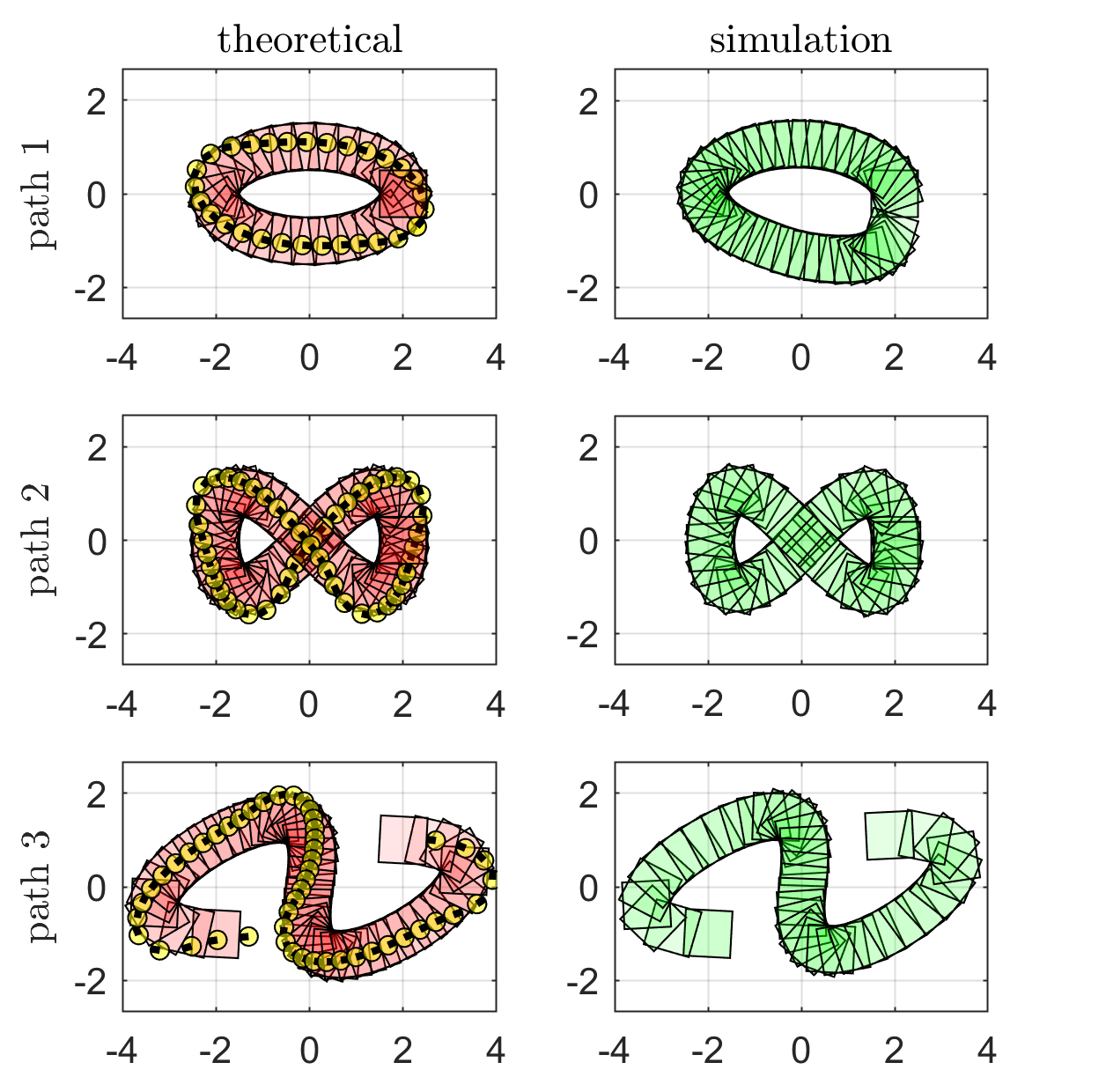}
\caption{Illustration of theoretical trajectories (left) and simulated trajectories (right) for the flat paths from Example \ref{ex:1}.}
\label{fig:flatPaths}
\end{figure}

\begin{figure}
\centering
\includegraphics[width=.5\columnwidth]{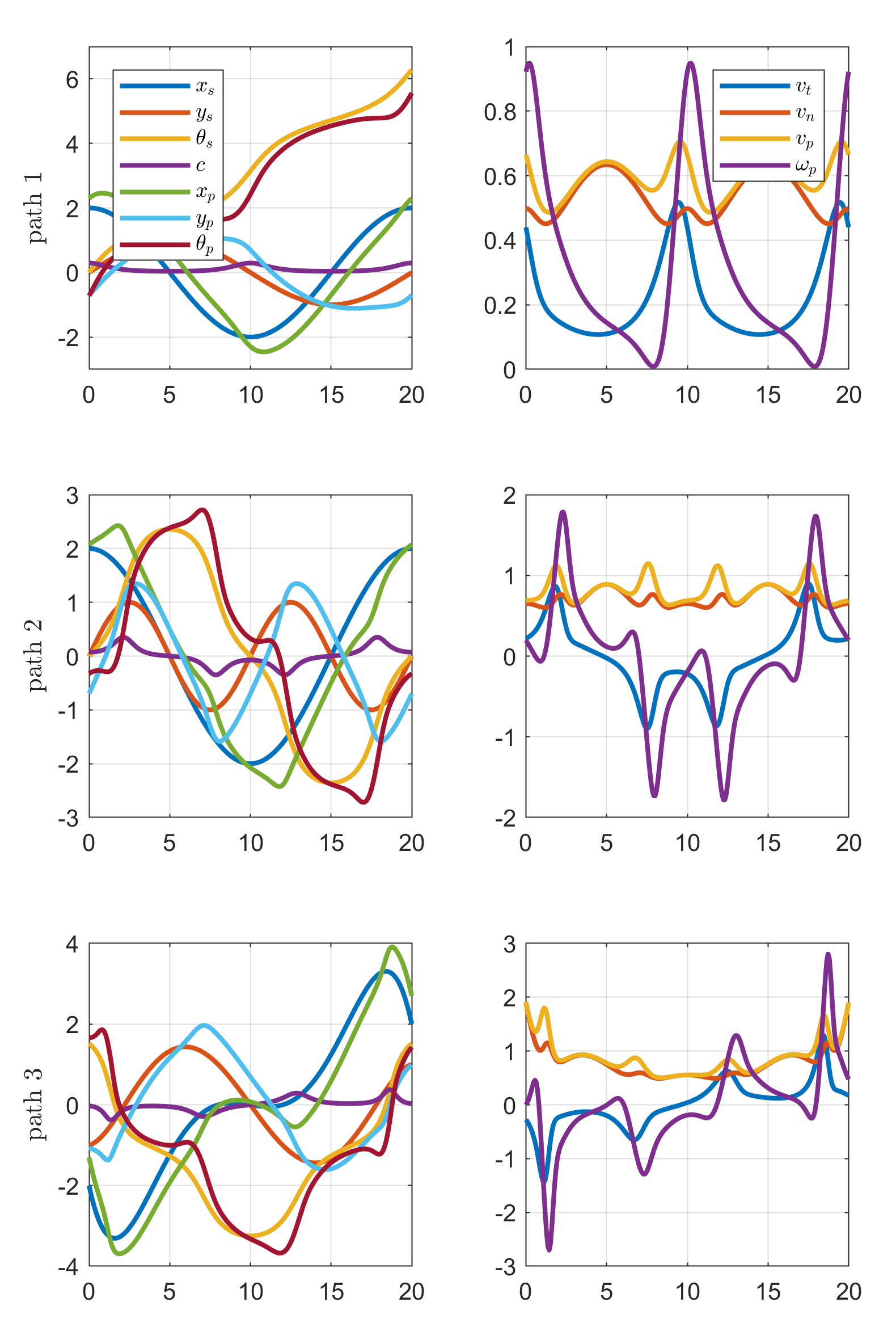}
\caption{Illustration of theoretical state trajectories (left) and input trajectories (right) for the flat paths from Example \ref{ex:1}.}
\label{fig:exp}
\end{figure}

\begin{example} 
\label{ex:1}
We give here a number of examples to demonstrate the state-action trajectories, $(\vectorstyle{x},\vectorstyle{u})(t)$, associated to several flat paths, $\vectorstyle{\zeta}(t)$. Three paths are illustrated in Fig. \ref{fig:flatPaths}. The associated state and action trajectories are given in Fig. \ref{fig:exp}. The third trajectory is constructed using interpolating B-splines of polynomial degree $5$ with 2 knots on the waypoint{s} given below (see section for details \ref{sec:geometric-path-planning}) \cite{stoical2016obstacle}. For all experiments, we set $T=20$ seconds. The geometrical parameters are given by $a=1$, $b=1$ and $r=2\cdot 10^{-1}$. 
\begin{equation*}
	\begin{aligned}
		\vectorstyle{\zeta}_1(t) &= \left(2\cos\left(\tfrac{2\pi}{T}t\right), \sin\left(\tfrac{2\pi}{T}t\right)\right), ~t \in[0,T] \\
		\vectorstyle{\zeta}_2(t) &= \left(2\cos\left(\tfrac{2\pi}{T}t\right),2\cos\left(\tfrac{2\pi}{T}t\right)\sin\left(\tfrac{2\pi}{T}t\right) \right), ~t \in[0,T]\\
		\vectorstyle{\zeta}_3(t_i) &\in \{(-2,-1),(-2,1),(0,1),(0,-1),(2,-1),(2,1)\}
	\end{aligned}
\end{equation*}

The expressions for the orientation and contact point can be verified to satisfy physical intuition. To validate the input expression we perform open-loop simulations. To stress the difference between the calculated and simulated trajectories we distort the input signals with Gaussian white noise ($\sigma=0.05$) on a time grid of $250$ intervals. The results are also depicted in Fig. \ref{fig:flatPaths}.
\end{example}

Theorem \ref{th:1} has value on its own provided the rich literature on flatness based control which now extends to slider-pushers. In the next section we discuss a method tailored to path planning that relies on an additional property of the system.


\section{Efficient Path Planning}
In this section we describe an original two step approach tailored to constrained time optimal collision free path planning. Our approach relies on a particular invariance property of the flat expressions for the slider-pusher system.

\subsection{Invariance to path parametrization}\label{sec:invariance-to-path-parametrization}

Here we demonstrate that the slider-pusher state is invariant to the time differential transformations of the path parametrizing coordinate. 

To explain what we mean by that, assume that a flat path, $\vectorstyle{\zeta}(\tau)$, as a function of some scalar path coordinate, $\tau$, has been obtained. The path coordinate, $\tau$, determines the spatial geometry of the \textit{flat path} but does not necessarily has to coincide with the time coordinate, $t$. The path’s actual time dependency follows from any relation, $\tau(t)$, between the path coordinate, $\tau$, and time coordinate, $t$. This relation determines the pace at which a system traverses the geometric path and how the system would be perceived correspondingly by a physical observer. 

Given that in the present context we are interested in the time derivatives of the path, we want to find a relation between the \textit{geometric} derivatives of the path to which we have access, and, the corresponding time derivatives. {For clarity we denote the geometric derivatives using prime notation, i.e. $\vectorstyle{\zeta}'$, $\vectorstyle{\zeta}''$, referring to $\frac{\text{d}}{\text{d}\tau}\zeta$, $\frac{\text{d}^2}{\text{d}\tau^2}\zeta$, $\dot{\vectorstyle{\zeta}}$, etc. and the time derivatives using dot convention $\ddot{\vectorstyle{\zeta}}$.} Relying on the chain rule, it is easily verified that the time derivatives of the path can be rewritten as a function of the geometric derivatives of the path and the time derivatives of the path coordinate, $\dot{\tau}$, $\ddot{\tau}$, etc. The relation between the time and geometric derivatives of the path is governed by a linear mapping. Though the linear map itself is a nonlinear function of $\dot{\tau}$, $\ddot{\tau}$, etc.
\begin{equation}
\label{eq:4}
\begin{aligned}
	\begin{pmatrix}
		\dot{\vectorstyle{\zeta}} \\ \ddot{\vectorstyle{\zeta}} \\ \dddot{\vectorstyle{\zeta}} \\ \vdots
	\end{pmatrix} &= \begin{pmatrix}
		\dot{\tau} & 0 & 0 & \cdots \\
		\ddot{\tau} & \dot{\tau}^2 & 0 &  \\
		\dddot{\tau} & 3 \dot{\tau} \ddot{\tau} & \dot{\tau}^3 & \\
		\vdots & &   & \ddots
	\end{pmatrix} \otimes \mathrm{I}_{n_\zeta} \cdot\begin{pmatrix}
		\vectorstyle{\zeta}' \\ \vectorstyle{\zeta} '' \\ \vectorstyle{\zeta}''' \\ \vdots
	\end{pmatrix} 
\end{aligned}
\end{equation}

Substituting this relation into the operator $\Phi$ suggests that the corresponding state trajectory is determined by the geometric derivatives of the path as well as the time derivatives of the relation between the path coordinate and the time coordinate itself. In general it can not be expected that the operator $\Phi$ will be invariant to the time derivatives of $\tau$. This means that if we were to change the pace at which we traverse the path, $\vectorstyle{\zeta}$, this would also influence the geometry of the state trajectory, $\vectorstyle{x}$.
\begin{equation}
\begin{aligned}
	\vectorstyle{x} &= \Phi(\vectorstyle{\zeta},\dot{\vectorstyle{\zeta}},\ddot{\vectorstyle{\zeta}},\dots) \\
	&= \Phi(\vectorstyle{\zeta},{\vectorstyle{\zeta}}',{\vectorstyle{\zeta}}'',\dots,\dot{\tau},\ddot{\tau},\dots)
\end{aligned}
\end{equation}

For the slider-pusher, $\vectorstyle{x}=(\vectorstyle{x}_s,\vectorstyle{x}_p)$, it can be verified by simple substitution that $\Phi$ is in fact invariant to $\tau$, {ergo
\begin{equation}
	\begin{aligned}
		\vectorstyle{x} = \Phi(\vectorstyle{\zeta},{\vectorstyle{\zeta}}',{\vectorstyle{\zeta}}'',\dots,\dot{\tau},\ddot{\tau},\dots) = \Phi(\vectorstyle{\zeta},{\vectorstyle{\zeta}}',{\vectorstyle{\zeta}}'',\dots) 
	\end{aligned}
\end{equation}}
{Henceforth} we will refer to $\dot{\tau},\ddot{\tau},\dots$ as the time differential properties of the path coordinate.{ The property is summarized in the following theorem.} {For completeness we refer to Appendix \ref{sec:verification-of-theorem-refth2} for an illustration.}
\begin{theorem}
\label{th:2}
The geometry of the state trajectories of the slider-pusher system {is} invariant to the time differential properties of the path coordinate, $\tau$.
\end{theorem}

This is not a general property of flat systems. 
\begin{example}
Consider the example of a gantry crane. For details about the model see \cite{fliess1995flatness}. The flat coordinate of the gantry crane is determined by the Cartesian position of the load while its full state also includes the angle between the load and the gantry. The relation between the angle and the flat coordinate is given by $\theta=\arctan \frac{\ddot{x}}{g-\ddot{y}}$. Substitution of $\ddot{x} = x'\ddot{\tau}+x''\dot{\tau}^2$ and $\ddot{y} = y'\ddot{\tau}+y''\dot{\tau}^2$ into the latter expression illustrates that the relative position of the load and gantry is determined by the time differential properties of the flat path. Hence we conclude that the flatness of the gantry crane is not invariant to any time differential transformation of the path parametrizing coordinate.
\end{example}

Without loss of generality, throughout we will further assume that the trajectory starts at $t = 0$, ends at $t = T$ and that $\tau(0) = 0 \leq \tau(t) \leq  \tau(T) = 1$.

\subsection{Imposing arbitrary velocity profiles}\label{sec:imposing-arbitrary-velocity-profiles} Clearly, the invariance of the slider-pusher state trajectory, $\vectorstyle{x}(t)$, does not extend to the input trajectory, $\vectorstyle{u}(t)$. Put differently, $\Psi$ is not invariant to the time differential properties of the path. Furthermore, although the state trajectory is invariant, the actual position of the system as a function of time still depends on the function $\tau$ implying that the relation between the path and time coordinates is still important. 

In this section we discuss how a meaningful choice for $\tau$ can be determined as a function of the velocity of the slider, $v_s$. Consider therefore the following differential relation between the throughput velocity, $v_s$, and the first time derivative of the path coordinate, $\tau$. The relation is governed by a factor {which we denote concisely using the variable $\psi$. Note that the variable $\psi$ is path dependent meaning it can be determined given $\vectorstyle{\zeta}$.} 

\begin{equation}
\label{eq:Vdtau}
\begin{aligned}
	v_s &= \sqrt{\dot{x}_s^2+\dot{y}_s^2} = \sqrt{(x_s')^2+(y_s')^2} \dot{\tau} = \frac{1}{\psi} \dot{\tau}  \\
	\dot{\tau} &= \psi v_s
\end{aligned}
\end{equation}

Higher order time derivatives of $\tau$ can be determined by deriving the expression above. Similar to the relation between the time and geometric derivatives of the path, the time derivatives of $\tau$ and $v_s$ are related through a linear mapping. The linear map itself is once more nonlinear in the geometric derivatives of the factor $\psi$. The values for $\dot{\tau}$, $\ddot{\tau}$, $\dots$, can be found by forward substitution.
\begin{equation}
\label{eq:tauV}
\begin{aligned}
	\begin{pmatrix}
		\dot{\tau} \\ \ddot{\tau} \\ \dddot{\tau} \\ \vdots
	\end{pmatrix} &= \begin{pmatrix}
		\psi & 0 & 0 & \cdots \\
		\psi'\dot{\tau} & \psi & 0 &  \\
		\psi''\dot{\tau}^2+\psi'\ddot{\tau} & 2 \psi'\dot{\tau} & \psi & \\
		\vdots & &   & \ddots
	\end{pmatrix} \cdot\begin{pmatrix}
		v_s \\ \dot{v}_s \\ \ddot{v}_s \\ \vdots
	\end{pmatrix} 
\end{aligned}
\end{equation}

These equations imply that we can impose an arbitrary velocity profile, $v_s(t)$, and find the corresponding relation $\tau$ by integration of $\psi v_s$ and evaluation of (\ref{eq:tauV}). Substituting these expressions, together with the path derivatives of $\vectorstyle{\zeta}$, into $\Psi$, will provide the associated input signals as a function of the geometric path and the time derivatives of the arbitrary velocity profile. Note that all the while the geometry of the state trajectory will not change.

\begin{figure*}[h!]
\centering
\includegraphics[width=\columnwidth]{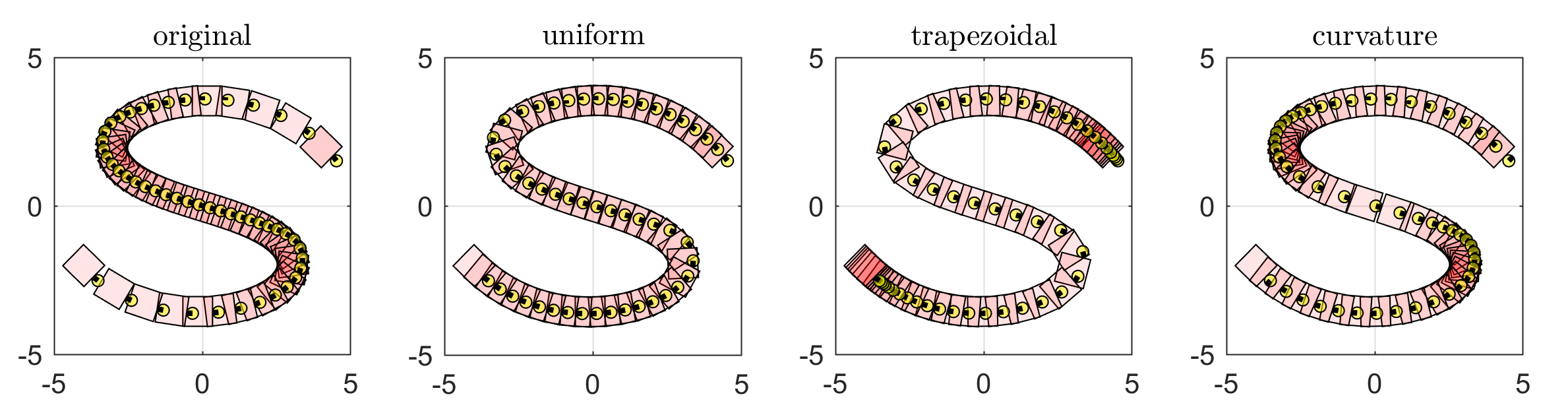}
\caption{Illustration of imposed velocity profiles from Examples \ref{ex:paths} and \ref{ex:paths2} on the flat path from Example \ref{ex:paths}. \textit{ From left to right:} original path, constant velocity, trapezoidal velocity and curvature dependent velocity.}
\label{fig:paths}
\end{figure*}

\begin{figure}
\centering
\includegraphics[width=.5\columnwidth]{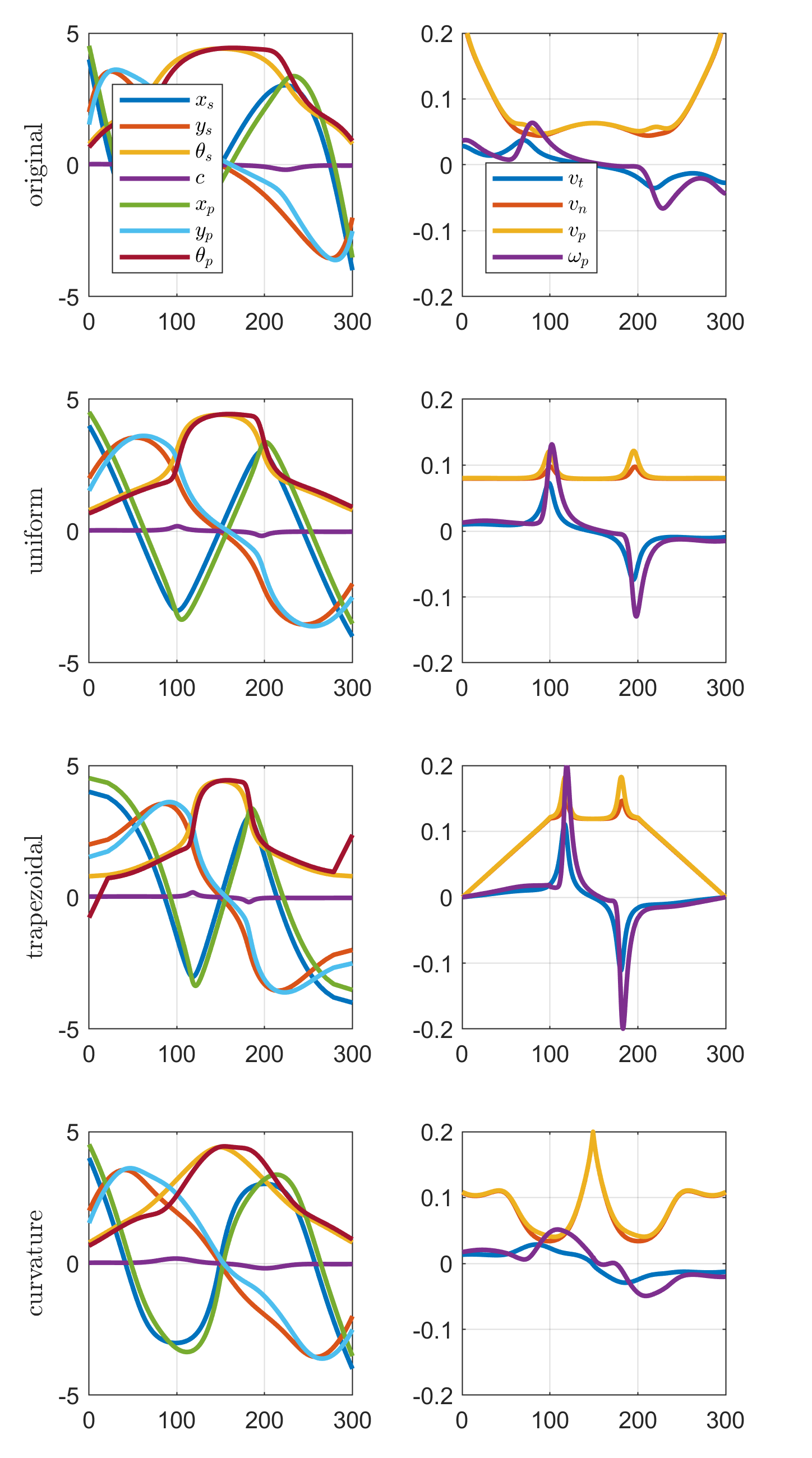}
\caption{Illustration of the state (left) and input (right) trajectories imposing the velocity profiles from Examples \ref{ex:paths} and \ref{ex:paths2} on the flat path from Example \ref{ex:paths}. Note how the state trajectories for the different velocity profiles are aliases warped over the horizontal axis whilst the input trajectories are clearly not.}
\label{fig:pathsV}
\end{figure}

\begin{example} 
\label{ex:paths}
As illustrative examples, one may consider uniform or trapezoidal velocity profiles for the slider.
\begin{equation*}
	\begin{aligned}
		v_s(t) &= v_0\\
		v_s(t) &= a_0 \left(\mathbf{1}_{[0,\Delta)}(t) t +  \mathbf{1}_{[\Delta,2\Delta)}(t)\Delta + \mathbf{1}_{[2\Delta,3\Delta]}(t)(3\Delta-t)\right)
	\end{aligned}
\end{equation*}

We use the arbitrary flat path
\begin{equation}
	\vectorstyle{\zeta}({\tau}) = \left\lbrace
	\begin{aligned}
		\left(4 \cos({\tau}),2 \sin({\tau}) + 2\right),&~ {\tau} \in \left[0,\frac{3\pi}{2}\right] \\
		\left(4 \cos({\tau}),2 \sin({\tau}) - 2\right),&~ {\tau} \in \left[\frac{3\pi}{2},3\pi\right] 
	\end{aligned}
	\right.
\end{equation}
for demonstration. 

In Fig. \ref{fig:paths}, the behavioural implications of imposing a certain velocity profile are visualized. The snapshots are taken equidistantly in time. One observes that although the system clearly visits different locations on the path at different time instances, the geometry of the path -- as reflected by the paths of the centre of mass of the slider and pusher and their orientation -- is unaffected by the velocity profile. This observation is confirmed by Fig. \ref{fig:pathsV}, where the state and input trajectories are represented. Clearly, the different state trajectories are aliases warped along the horizontal axis, whilst the different input trajectories are not.
\end{example}

\begin{example}
\label{ex:paths2}
Here we impose a proportional relation between a geometrical property of the path, such as its curvature, $\kappa$, and the slider velocity, $V$. In this case there is no velocity profile available since the profile depends on the geometric properties of the path. 

For illustrational purposes we propose the following rule where the hyper parameter $\kappa_0$ prevents the denominator from going to zero and hence the velocity to infinity.
\begin{equation*}
	\frac{1}{\kappa_0+\kappa} \propto V = \frac{1}{\psi} \dot{\tau}
\end{equation*}

Now, the curvature, $\kappa$, can be expressed as a function of $\vectorstyle{\zeta}$
\begin{equation*}
	\kappa = \sqrt{\left(\frac{\text{d}^2 x_s}{\text{d} s^2}\right)^2+\left(\frac{\text{d}^2 y_s}{\text{d} s^2}\right)^2}
\end{equation*}
where $s$ represents the arc length. 

To complete the analysis we now require a relation between the arc length, $s$, and the path coordinate, $\tau$. The latter can be obtained in analogy to the reasoning behind (\ref{eq:Vdtau}) and (\ref{eq:tauV}).
\begin{equation}
	\begin{aligned}
		\frac{\text{d}\tau}{\text{d}s} &= \psi \frac{\text{d}s}{\text{d}s} = \psi\\
		\frac{\text{d}^2\tau}{\text{d}s^2} &= \psi' \frac{\text{d}\tau}{\text{d}s} = \psi'\psi
	\end{aligned}
\end{equation}
so that 
\begin{equation*}
	\dot{\tau} = \frac{\psi }{\kappa_0 + \sqrt{\left(\psi^2 X'' + \psi'\psi X'\right)^2 + \left(\psi^2 Y'' + \psi'\psi Y'\right)^2 }}  = \chi
\end{equation*}

In Fig. \ref{fig:paths} and Fig. \ref{fig:pathsV} the analysis from example \ref{ex:paths} is repeated. As a result the velocity of the slider decreases in corners and accelerates when the motion becomes rectilinear again. 
\end{example}

In order to impose the condition $\tau(T)= 1$, or equivalently $t(1) = T$ for arbitrary $T$, we may substitute $\eta \hat{v}_s$ for the velocity profile $v_s$ where $\hat{v}_s(t)$ now determines a velocity primitive and $\eta$ is the required scaling factor. This strategy can be used to impose $T$ regardless of any of the methods described so far. 

From (\ref{eq:tauV}) it follows that $\text{d}\tau = \eta \hat{v}_s(t)\text{d}t$ so that
\begin{equation}
\eta \cdot \int_{0}^{T} \hat{v}_s(t)\text{d}t =\int_{0}^{1}\frac{1}{\psi(\tau)}\text{d}\tau
\end{equation}

The invariance property and the associated strategy to impose arbitrary velocity profiles are advantageous features of the slider-pusher system with particularly interesting implications for certain path planning problems. It follows that we can use some arbitrary path planning strategy to compute a reference trajectory for the slider-pusher system. Once this path is available, an auxiliary method can be applied to impose an arbitrary velocity profile, or, determine a path dependent velocity profile and thus meets some additional requirements. Such an approach is discussed next.

\subsection{Constrained time optimal collision free path planning}

We will argue here that the invariance property discussed in section \ref{sec:invariance-to-path-parametrization}, allows to decompose a constrained time optimal collision free path planning problem, approximately, into two easier and independent subproblems. 

The problem that we aim to solve is the following
\begin{equation}
\label{eq:ctocfpp}
\begin{aligned}
	\min_{\vectorstyle{x},\vectorstyle{u}} T \text{ s.t. } 
	0&=\vectorstyle{f}(\vectorstyle{x}(t),\vectorstyle{u}(t))-\dot{\vectorstyle{x}}(t)  \\
	0&\leq \vectorstyle{g}(\vectorstyle{x}(t)) \\
	0&\leq \vectorstyle{h}(\vectorstyle{u}(t)) \\
	{0} & {=\vectorstyle{k}(\vectorstyle{x}(0),\vectorstyle{x}(T))\color{red}} 
\end{aligned}
\end{equation}
{where $\vectorstyle{x}$ and $\vectorstyle{u}$ denote the state and input of the state-space representation of the slider-pusher. The vector functions, $\vectorstyle{g}$, $\vectorstyle{h}$, and $\vectorstyle{k}$, determine path constraints and boundary conditions respectively.} The dynamic constraint restricts the solution to the feasible state-action function space, $\mathcal{F}$, also see (\ref{eq:feasible}). For our arguments to apply, it is important that the constraints on the state and input trajectories are mutually independent and that the path constraints are geometric, i.e. do not depend on the particular time instant, expect for the boundary conditions at $t=0$ and $t=T$.

Now recall that we are interested in path planning problems for the slider-pusher. Hence we may exploit the differential flatness of the system to get rid of the differential constraint and reparametrize optimization problem (\ref{eq:ctocfpp}) using a flat trajectory.
\begin{equation}
\label{eq:ctocfpp2}
\begin{aligned}
	\min_{\vectorstyle{\zeta}} T \text{ s.t. } 
	0&\leq \vectorstyle{g}(\Phi(\vectorstyle{\zeta},\dot{\vectorstyle{\zeta}},\ddot{\vectorstyle{\zeta}},\dots) ) \\
	0&\leq \vectorstyle{h}(\Psi(\vectorstyle{\zeta},\dot{\vectorstyle{\zeta}},\ddot{\vectorstyle{\zeta}},\dots) )
\end{aligned}
\end{equation}

Second, we consider the following reformulation of $T$
\begin{equation}
\label{eq:T}
T = \int\nolimits_0^T \text{d}t = \int\nolimits_0^1 \frac{1}{\dot{\tau}}\text{d}\tau 
\end{equation}

Finally, we substitute $T$ into problem (\ref{eq:ctocfpp2}) and consider an arbitrary coordinate parametrization of the flat path. Then we obtain the following optimization problem parametrized by $\vectorstyle{\zeta}$ and ${\tau}$. Note that here we used the invariance property of the slider-pusher state.
\begin{equation}
\label{eq:ctocfpp3}
\begin{aligned}
	\min_{\vectorstyle{\zeta},\tau} \int\nolimits_0^1 \frac{1}{\dot{\tau}}\text{d}\tau \text{ s.t. } 
	0&\leq \vectorstyle{g}(\Phi(\vectorstyle{\zeta},{\vectorstyle{\zeta}}',{\vectorstyle{\zeta}}'',\dots) ) \\
	0&\leq \vectorstyle{h}(\Psi(\vectorstyle{\zeta},{\vectorstyle{\zeta}}',{\vectorstyle{\zeta}}'',\dots,\dot{\tau},\ddot{\tau},\dots) )
\end{aligned}
\end{equation}

To solve this problem efficiently, we will invoke the following approximation. According to (\ref{eq:tauV}) we have that $\dot{\tau} = \psi v_s$. If we then further assume that $v_s$ is constant, the objective, $T$, becomes independent of time and hence independent of the relation between $\tau$ and $t$.

As such we obtain an optimization problem where the objective, say $f(x)$, and one set of constraints, say $0 \leq g(x)$ depend on a certain set of variables, say $x$, and, an auxiliary set of constraints, say $0\leq h(x,y)$, that depend on the first set of variables, $x$, but also on an second set of variables, say $y$. It is important to note that the objective itself does not depend on the secondary set of variables, $y$. {Note that in the present setting $x$ refers to the geometric path $\zeta$ whilst $y$ refers to $\tau$ as a function of $t$.}

Under the right conditions we may {then} solve the optimization problem without the second set of constraints for the first set of variables. The second set of variables can then be found by finding values for the second set of variables so that the second set of constraints are satisfied. Equivalently, we have $\min_{x,y} f(x) \text{ s.t. } 0\leq g(x) \wedge 0\leq h(x,y)$. We solve for $x$ first, yielding $x^*$, whilst neglecting constraint, $h$. Then we find some $y$ so that $0\leq h(x^*,y)$. {Provided that $h$ is sufficiently well behaved, there will exist $y^*$ so that $0\leq h(x^*,y^*)$. For $\dim h(x^*,\cdot) = \dim y$, $h$ must simply be invertible. For $\dim h(x^*,\cdot) < \dim y $ the problem is undetermined and a secondary objective can be included.}

This analysis suggests that we can solve problem (\ref{eq:ctocfpp3}) in two steps. In the first step, we solve for $\vectorstyle{\zeta}$, omitting the second set of constraints. In the second step, we solve for $\tau$, correcting the approximation made above, whilst recycling the solution for $\vectorstyle{\zeta}$ from the first step. This approach yields a strictly geometric problem and strictly time dependent problem.

\begin{figure*}[h!]
\centering
\includegraphics[width=\columnwidth]{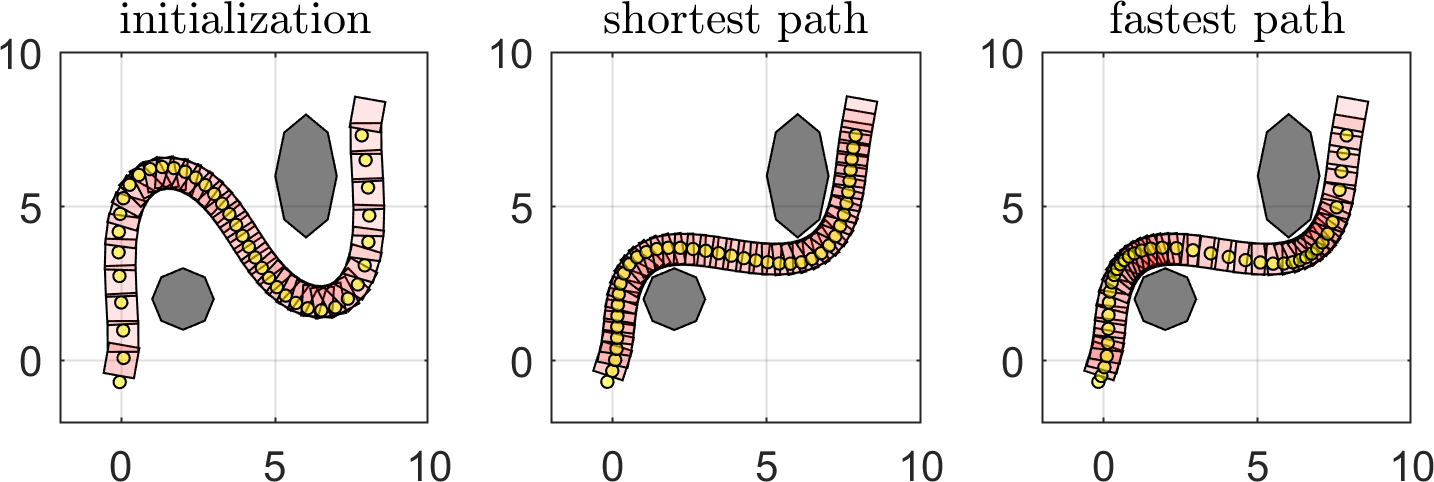}
\caption{Illustration of the two step constrained time optimal collision free path planning approach. We initialize the solution with some arbitrary collision free path parametrized by the path coordinate $\tau$ (left). Then we solve the geometric path planning problem from section \ref{sec:geometric-path-planning}, generating the shortest path (center). Second, we solve the constrained time optimization problem from section \ref{sec:constrained-time-optimization}, generating the fastest path that satisfies the input constraints (right). Note how the slider decelerates in the advent of a corner to respect $\omega_p<\overline{\omega}_p$ and `catches up' by accelerating on straight strips, up to $v_p = \overline{v}_p$.}
\label{fig:flatPathOpt}
\end{figure*}

\subsubsection{Geometric path planning}\label{sec:geometric-path-planning} For the first subproblem we aim to solve the following geometric path planning problem which computes the shortest collision-free path. It is important to remark that the objectives and constraints depend on the geometry of the path alone.

Therefore the interpretation of the path coordinate is irrelevant and we can solve for arbitrary $\tau$. 
\begin{equation}
\label{eq:opt1}
\vectorstyle{\zeta}^* = \arg \min_{\vectorstyle{\zeta}} \int\nolimits_{0}^1 \frac{1}{\psi} \text{d}\tau \text{ s.t. } 0\leq \vectorstyle{g}(\Phi(\vectorstyle{\zeta},{\vectorstyle{\zeta}}',{\vectorstyle{\zeta}}'',\dots) )
\end{equation}

Now remember that $\psi^{-1} = \sqrt{(x_s')^2+(y'_s)^2}$. It follows that this problem strictly reads as minimizing the path length whilst satisfying the auxiliary path constraints. This in fact is unsurprising, since by definition, the fastest path is in fact the shortest path considering that we travel at constant speed.

\subsubsection{Constrained time optimization}\label{sec:constrained-time-optimization} For the second subproblem we aim to optimize the time required for the system to execute the geometrical path obtained by solving (\ref{eq:opt1}) whilst now also satisfying the input constraints. Though instead of finding just some $\tau$ that satisfies the input constraints, we will correct for the approximation made whilst addressing problem (\ref{eq:ctocfpp3}). Our objective is now to minimize (\ref{eq:T}) whilst satisfying input constraints for $\vectorstyle{\zeta}^*$.

Inspired by the work of \cite{debrouwere2013convex,verscheure2009time}, we recast this as a standard optimal control problem. Therefore we introduce the auxiliary state variable $\vectorstyle{z}=(z,z',z'')$ where  $z=\dot{\tau}^2$ and the auxiliary control variable $v = z'''$. Again the accent notation denotes derivatives to the path coordinate. It follows that the time derivatives of $\tau$ can be expressed as a nonlinear function of $\vectorstyle{z}$ and $v$. One verifies
\begin{equation}
\begin{aligned}
	\dot{\tau} &= \sqrt{z}\\
	\ddot{\tau} &= \tfrac{1}{2}z'\\
	\dddot{\tau} &= \tfrac{1}{2}\sqrt{z}z'' \\
	\ddddot{\tau} &= \tfrac{1}{2}zv+\tfrac{1}{4}z'' 
\end{aligned}
\end{equation}

Based on these auxiliary variables, we can rewrite the objective (\ref{eq:T}) as a function of the state variable, $z$, and the input constraint as a function of the state $\vectorstyle{z}$ and the control $v$
\begin{equation}
\label{eq:opt2}
\begin{aligned}
	\min_{\vectorstyle{z},v} \int_{0}^1 \frac{1}{\sqrt{z}}\text{d}\tau \text{ s.t. }&\vectorstyle{z}' = \matrixstyle{A}\vectorstyle{z} + \matrixstyle{B}v\\
	0&\leq \vectorstyle{h}(\Psi(\vectorstyle{\zeta},{\vectorstyle{\zeta}}',{\vectorstyle{\zeta}}'',\dots,z,z',z'',v))
\end{aligned}
\end{equation}
which can be solved using standard solution techniques.

In conclusion, we remark that the relation between $z$ and $\dot{\tau}$ is arbitrary. This particular relation is motivated by what happens when either $\dot{\tau}$ or $z$ would approach zero. Therefore we want $z$ and $\dot{\tau}$ to be proportional and not inversely proportional. The squared relation is motivated by the objective. Consider therefore the primitive $F(x)$ of $f(x)^{-1}$. The simplest or most well-behaved function for which $f(x^*) = 0$ and $F(x^*)=0$ is $f(x)=\sqrt{x}$. This prevents the objective in (\ref{eq:opt2}) from being ill-defined when $z(\tau) = 0$.

\begin{figure}[h!]
\centering
\includegraphics[width=.5\columnwidth]{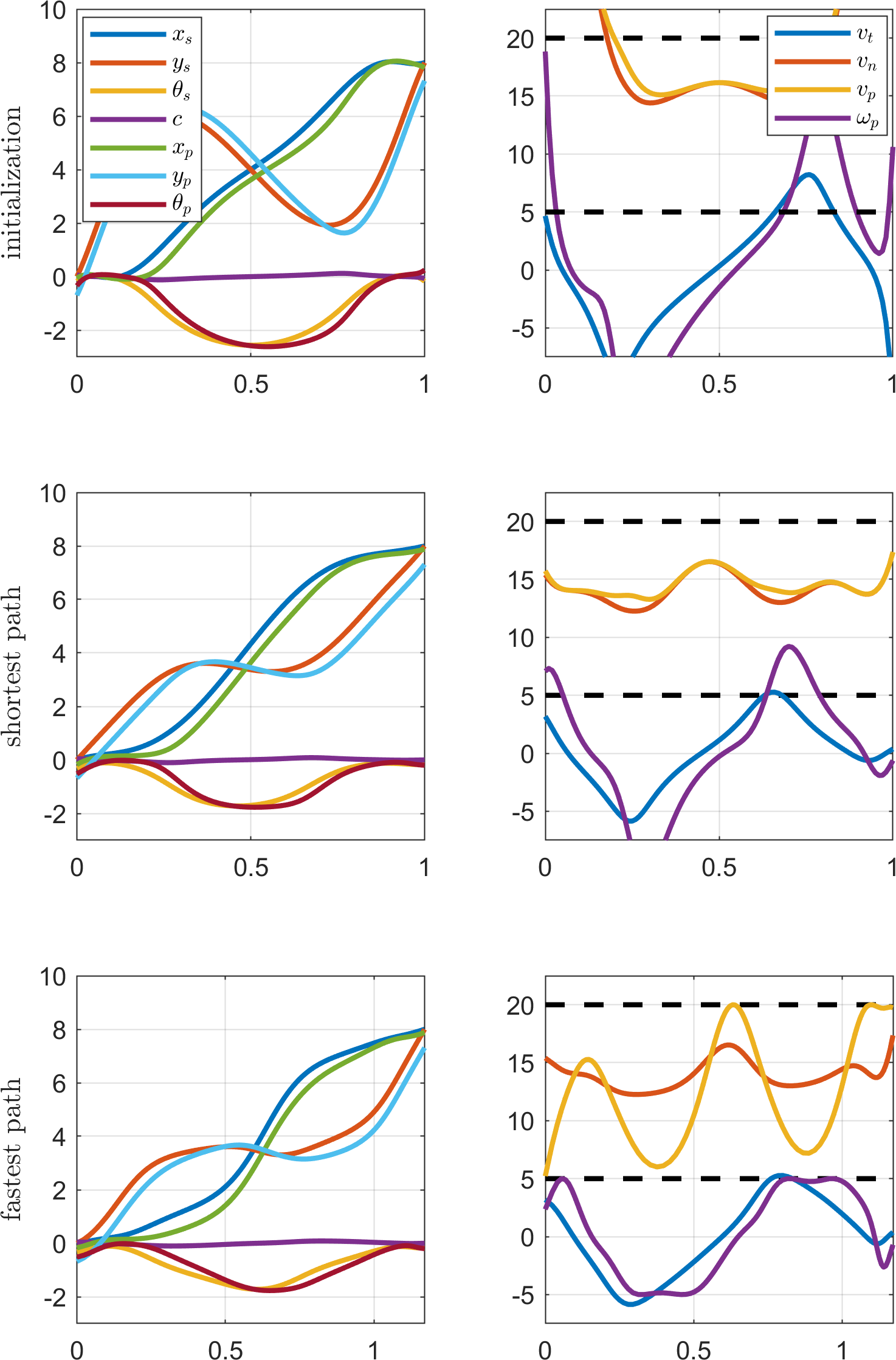}
\caption{Visualization of the state and input trajectories corresponding the paths in Fig. \ref{fig:flatPathOpt}. Remark how the shortest path violates the input constraints as opposed to the velocity constrained fastest path similar to the results depicted in Fig. \ref{fig:paths}.}
\label{fig:flatPathOptV}
\end{figure}

\subsection{Experiments}
Consider the path planning problem
\begin{equation}
\begin{aligned}
	\min_{\zeta,\tau} T \text{ s.t. }  0&= \vectorstyle{\zeta}(0) - (0,0)\\ 
	0&= \vectorstyle{\zeta}(1) - (8,8)  \\
	10^{-2} &\leq \mathrm{dist}(\mathrm{slider}({\Phi(\vectorstyle{\zeta},\vectorstyle{\zeta}',\vectorstyle{\zeta}'',\dots)}),\mathrm{obs})  \\
	{\vectorstyle{1}_c^\top \Phi(\vectorstyle{\zeta},\vectorstyle{\zeta}',\vectorstyle{\zeta}'',\dots) } & {\leq \tfrac{a}{2}} \\
	\underline{\vectorstyle{u}} &\leq {\Psi(\vectorstyle{\zeta},\vectorstyle{\zeta}',\vectorstyle{\zeta}'',\dots,\dot{\tau},\ddot{\tau},\dots)} \leq \overline{\vectorstyle{u}}
\end{aligned}
\end{equation}
Here, the function $\mathrm{dist}(\cdot,\cdot)$ compute the minimal distance between two polygons, the function $\mathrm{slider}(\cdot)$ computes a polygonal representation of the slider as a function of the full state and the object $\mathrm{obs}$ contains a polygonal representation of any obstacles. {The fourth constraint assures that the pusher never moves past a corner of the rectangle which would imply a loss of contact.} For the purpose of demonstration, we only consider constraints for the pusher motion, $|v_p|<20$ and $|\omega_p|<5$,{ and for the slider motion, $v_n \geq 0$}. The particularities of the polygon obstacles are given in the left frame of Fig. \ref{fig:flatPathOpt}. {Further note that here we could also include constraints affecting the full geometric state of the system, namely $\Phi(\vectorstyle{\zeta}(0),\vectorstyle{\zeta}'(0),\vectorstyle{\zeta}''(0),\dots)$ and $\Phi(\vectorstyle{\zeta}(1),\vectorstyle{\zeta}'(1),\vectorstyle{\zeta}''(1),\dots)$.}

We implemented the optimization problems (\ref{eq:opt1}) and (\ref{eq:opt2}) in \texttt{Matlab}. Both were solved with the SQP solver from FMINCON. Additional numerical details are shared for each subproblem before discussing results.

\subsubsection{Geometric path planning}\label{sec:geometric-path-planning2}
To solve problem (\ref{eq:opt1}) numerically, first we parametrize an arbitrary flat path, $\vectorstyle{\zeta}$, using B-splines. The number of splines, $N=m+d-1$, is determined by the number of knots, $m$, and the degree, $d$. In our experiments we use $(m,d) = (5,5)$. For details on B-splines in the context of flatness based path planning methods we refer to \cite{stoical2016obstacle,stoican2017constrained}. 
\begin{equation}
\vectorstyle{\zeta}(\tau;\vectorstyle{\theta}) = \sum\nolimits_{i=1}^N  \mathrm{B}_{i,d}(\tau) \vectorstyle{\theta}_i
\end{equation}

The integral is evaluated using trapezoidal integration over $10^2$ intervals.  We found empirically that it is useful to add $\gamma \psi^{-2}$ to the objective with $\gamma=10^{-1}$ for regularization. The collision constraints are imposed on the corresponding path coordinate grid. 

\subsubsection{Constrained time optimization}
Problem (\ref{eq:opt2}) can be treated using standard optimal control solution techniques. In particular we use a direct transcription method introducing the variables $\{\vectorstyle{z}_k,v_k\}_{k=0}^{K-1}$ over an equidistant grid $\{\tau_k\}_{k=0}^{K-1}$. 

Since the auxiliary system dynamics are linear, these can be evaluated exactly over the discrete grid. The input constraints are evaluated on the grid only. For numerical stability we do not approximate the integral using a trapezoidal rule. Instead we approximate $z(\tau)$ using a piecewise linear function between nodes $\tau_k$ and $\tau_{k+1}$. That way $\dot{\tau}_0 = z_0=\dot{\tau}_K = z_K=0$.
\begin{equation}
\label{eq:approx}
\int_{\tau_k}^{\tau_{k+1}} \frac{1}{\sqrt{z}(\tau)} \text{d}\tau \approx 2\frac{\tau_{k+1}-\tau_k}{\sqrt{z_{k+1}}+\sqrt{z_k}}
\end{equation}

\subsubsection{Results}
\begin{figure*}[h!]
\centering
\begin{subfigure}[b]{0.49\textwidth}
	\includegraphics[width=\textwidth]{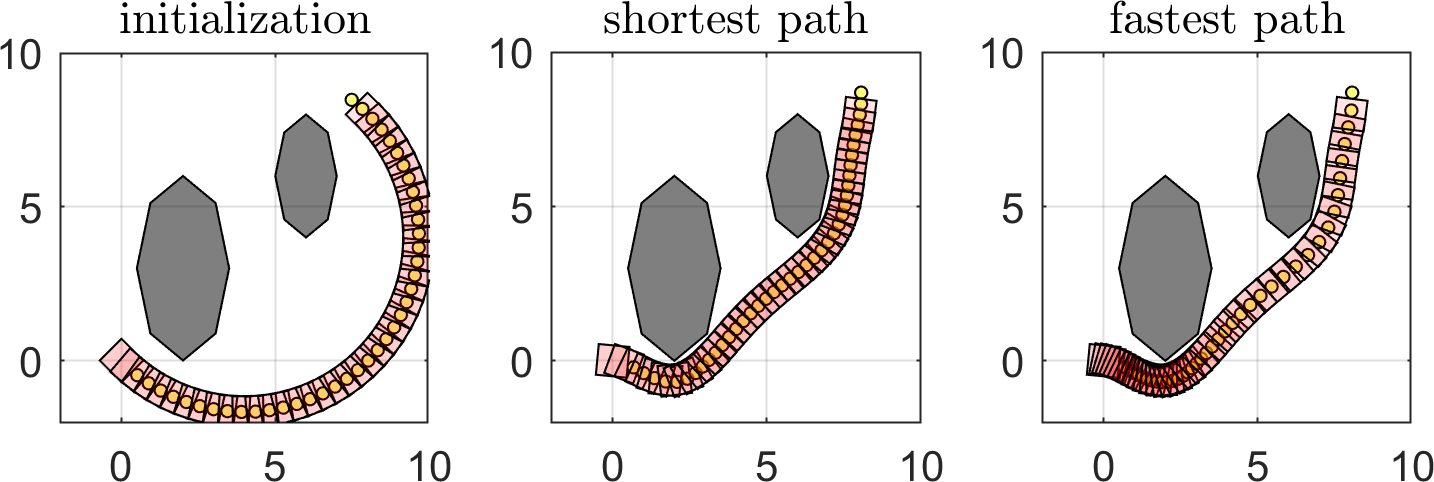}
	\caption{Second experiment}
\end{subfigure}
\begin{subfigure}[b]{0.49\textwidth}
	\includegraphics[width=\textwidth]{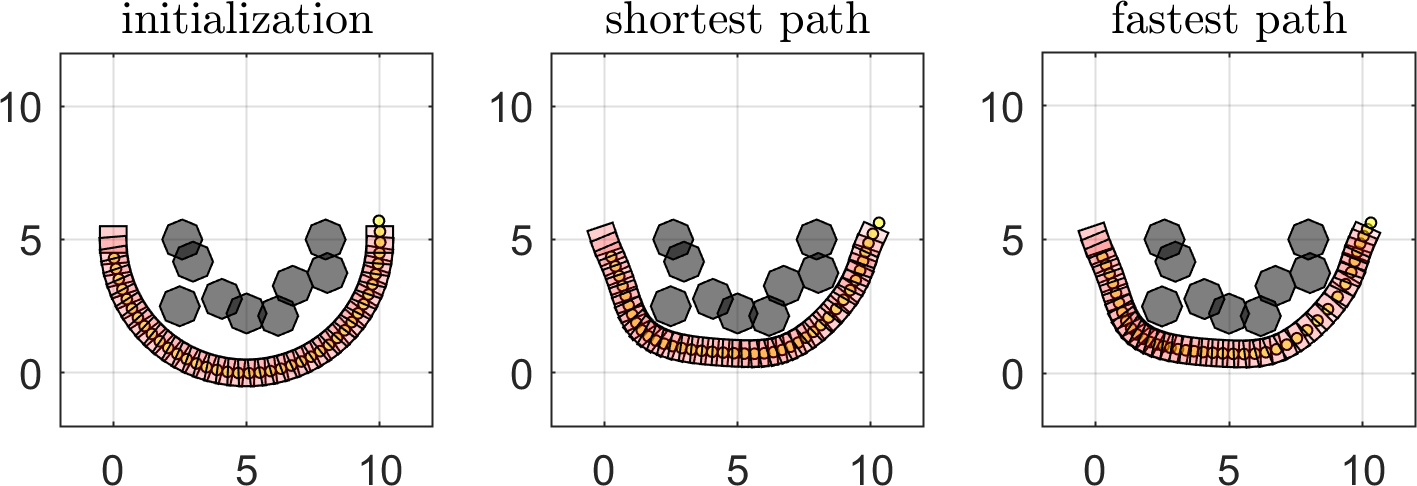}
	\caption{Third experiment}
\end{subfigure}
\caption{Illustration of additional experiments for the two step constrained time optimal collision free path planning approach.}
\label{fig:exp_paths}
\end{figure*}

\begin{figure*}[h!]
\centering
\begin{subfigure}[b]{0.49\textwidth}
	\includegraphics[width=\textwidth]{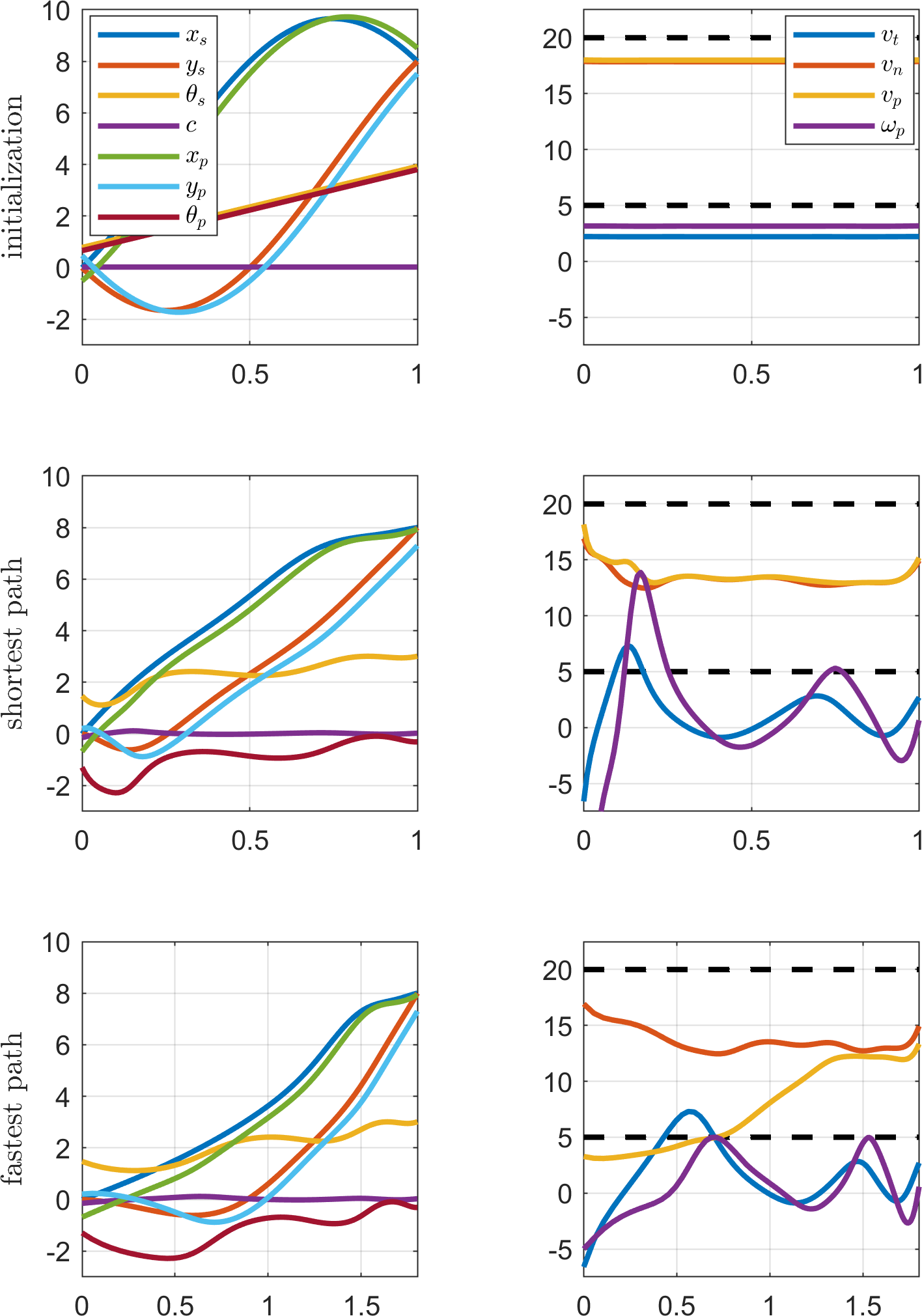}
	\caption{Second experiment}
\end{subfigure}
\begin{subfigure}[b]{0.49\textwidth}
	\includegraphics[width=\textwidth]{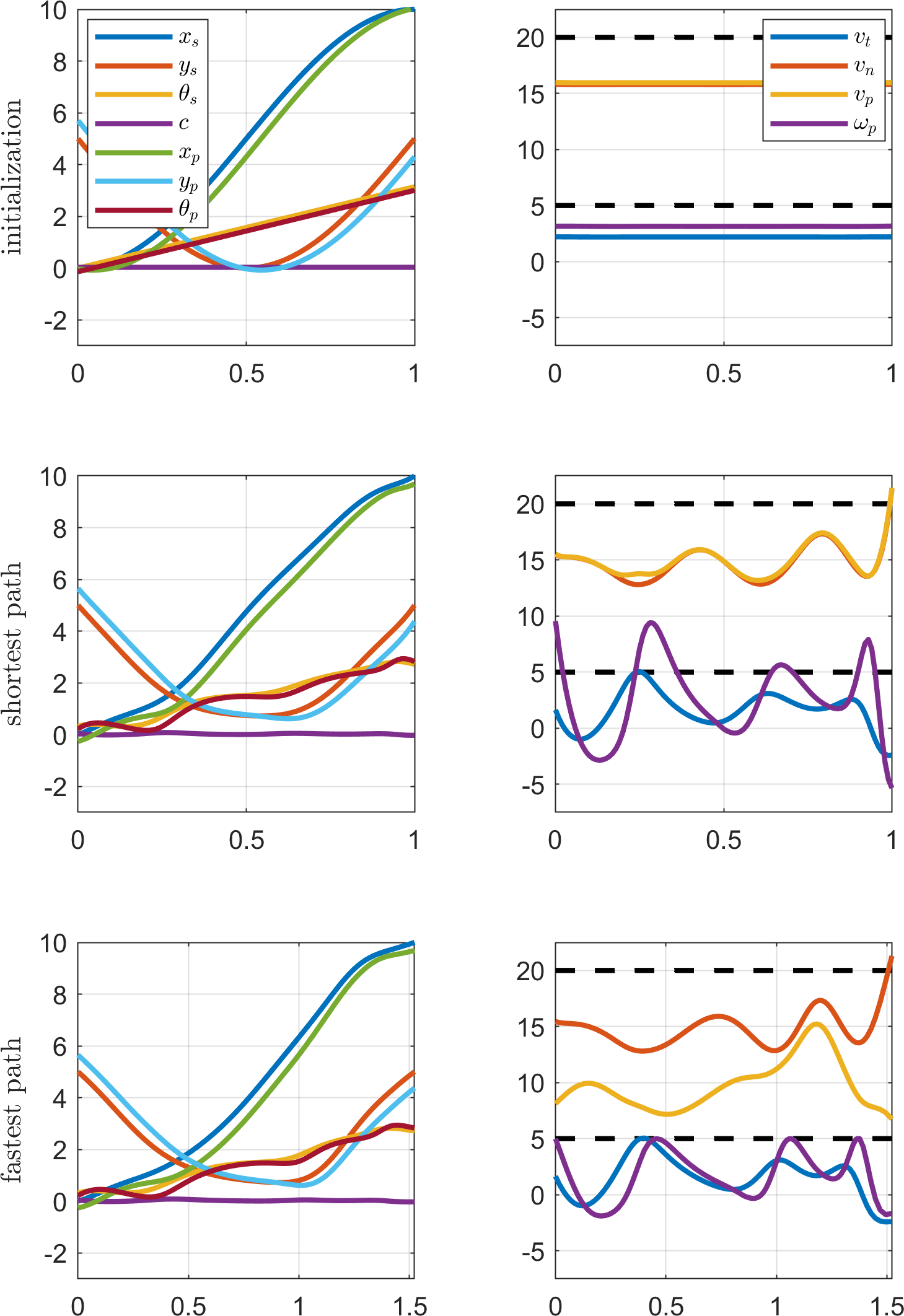}
	\caption{Third experiment}
\end{subfigure}
\caption{Illustration of additional experiments for the two step constrained time optimal collision free path planning approach.}
\label{fig:exp_velocities}
\end{figure*}

Results are visualized in Fig. \ref{fig:flatPathOpt} and Fig. \ref{fig:flatPathOptV}. Additional experiments are included in Fig. \ref{fig:exp_paths} and \ref{fig:exp_velocities}. The latter experiments are discussed in lesser detail.

Fig. \ref{fig:flatPathOpt} visualizes the initial path, shortest path (sec. \ref{sec:geometric-path-planning}) and fastest path (sec. \ref{sec:constrained-time-optimization}). Note how the slider decelerates in the advent of a corner to respect $\omega_p<\overline{\omega}_p$ and `catches up' by accelerating on straight strips, up to $v_p = \overline{v}_p$. Further, it can be verified that the influence of the regularization, as discussed in sec. \ref{sec:geometric-path-planning2}, is negligible and the numerical implementation of the geometric optimization problem retrieves the shortest path. 

The reciprocal dimensionless path coordinate is still normalized. Interpreting the path coordinate as the time coordinate would clearly result in a violation of the velocity constraints, even though improvement is present with respect to the initial trajectory. This is clearly visible in Fig. \ref{fig:flatPathOptV}.

In conclusion it can then be verified that the a posterior time optimization successfully satisfies the constraints whilst only slightly increasing the total time $T$ to $1.1699$.

\section{Conclusion}
In this contribution we have demonstrated that the quasi-static model of slider-pusher systems is differentially flat when frictionless contact between the slider and pusher is assumed {and the slider geometry is rectangular}. This has important implications with respect to various control tasks that benefit from flatness, specifically path planning problems. 

We have shown that the slider-pusher system in particular also satisfies an invariance property towards the time derivatives of the path coordinate. This invariance property has as a result that arbitrary velocity profiles can be imposed on any geometric path without suddenly violating any auxiliary geometric constraints on the state. This property was used to decompose a constrained time optimal path planning problems into two easier subproblems.

It is anticipated that these results will allow for improved and accelerated control design for slider-pusher systems. An interesting outlook is to exploit these properties in a model based predictive real-time control architecture {or in flatness based tracking control architectures are as well studied for other flat systems.
}

\appendix

{	\section{Generalised slider geometry}\label{sec:generalised-slider-geometry}
Here we briefly treat the case of a generalised slider geometry. As described in section \ref{sec:quasi-static-model} a slider with a smooth circumference can be parametrized using the angle $\phi$ and a distance function $r(\phi)$. The parametrization is so that the centre of mass is characterised by zero distance. Equivalently we have that
\begin{equation}
	\begin{aligned}
		\int_0^{2\pi} \int_0^{r(\phi)} r^2 \cos\phi \text{d}r\text{d}\phi &= 0 \\
		\int_0^{2\pi} \int_0^{r(\phi)} r^2 \sin\phi \text{d}r\text{d}\phi &= 0 \\
	\end{aligned}
\end{equation}

As a consequence of the properties of this parametrization it follows that objective in (\ref{eq:objective}) again generalizes to
\begin{equation}
	J \propto (\dot{x}_s^2+\dot{y}_s^2) + \beta_r^2 \theta_s^2
\end{equation}
where
\begin{equation}
	\beta_r^2  = \int_0^{2\pi} \int_0^{r(\phi)} r^3 \text{d}r\text{d}\phi
\end{equation}

In this case the contact point can be parametrized by a contact angle $\phi$. The corresponding contact distance is then given by $r(\phi)$ and the local position coordinate of the contact by
\begin{equation}
	\vectorstyle{p}_c^l = \begin{pmatrix}
		s_\phi r_\phi \\
		-c_\phi r_\phi 
	\end{pmatrix}
\end{equation}
so that
\begin{equation}
	\dot{\vectorstyle{p}}_c^l = \begin{pmatrix}
		c_\phi r_\phi \dot{\phi}+ s_\phi r_\phi'\dot{\phi}\\
		s_\phi r_\phi\dot{\phi}-c_\phi r_\phi' \dot{\phi}
	\end{pmatrix}
\end{equation}

We can then apply the same modelling strategy as described in section \ref{sec:principle-of-least-work}. Though first we give some additional attention to the contact point. In particular we define a local frame of reference that is locally normal to the slider's circumference at the contact point. This local frame of reference is characterised through three subsequent rotation's. The first rotation is the global orientation of the slider, $\theta$. The second rotation is determined by the slider circumference parametrization, $\phi$, in this case parametrizing the contact point. The third rotation is a result of the local nonlinearity of the circumference. As a result of the nonlinear function $r(\phi)$, an additional rotation $\alpha$ should be applied. This additional rotation is characterized as follows
\begin{equation}
	\tan \alpha = -\frac{r_\phi'}{r_\phi}
\end{equation}

Correspondingly
\begin{equation}
	\begin{aligned}
		\cos \alpha &= \frac{r_\phi}{\sqrt{r_\phi^2+ (r_\phi')^2}} \\
		\sin \alpha &= -\frac{r_\phi'}{\sqrt{r_\phi^2+ (r_\phi')^2}}
	\end{aligned}
\end{equation}

If we describe the inputs $\vectorstyle{u}$ in the frame of reference described above, i.e. 
\begin{equation}
	\begin{aligned}
		\vectorstyle{u} &= \matrixstyle{R}_\phi \tilde{\vectorstyle{u}} \\
		\tilde{\vectorstyle{u}} &= \matrixstyle{R}_\alpha \tilde{\tilde{\vectorstyle{u}}}
	\end{aligned}
\end{equation}
instead of the frame of reference attached to the slider itself, we obtain the following differential model.
\begin{equation}
	\begin{aligned}
		\dot{x} &= -\frac{\beta_r^2 r_\phi^2 + \beta_r^2 (r_\phi')^2}{\beta_r^2 r_\phi^2 + \beta_r^2 (r_\phi')^2 + r_\phi^2 (r_\phi')^2} s_{\alpha+\phi+\theta} \tilde{\tilde{v}}_n \\
		\dot{y} &= \frac{\beta_r^2 r_\phi^2 + \beta_r^2 (r_\phi')^2}{\beta_r^2 r_\phi^2 + \beta_r^2 (r_\phi')^2 + r_\phi^2 (r_\phi')^2} c_{\alpha+\phi+\theta} \tilde{\tilde{v}}_n \\
		\dot{\theta} &= \frac{r_\phi r_\phi'\sqrt{r_\phi^2+ (r_\phi')^2}}{\beta_r^2 r_\phi^2 + \beta_r^2 (r_\phi')^2 + r_\phi^2 (r_\phi')^2}\tilde{\tilde{v}}_n \\
		\dot{\phi} &= \frac{1}{\sqrt{r_\phi^2+ (r_\phi')^2}} \tilde{\tilde{v}}_t - \frac{r_\phi^3 r_\phi'}{\sqrt{r_\phi^2+ (r_\phi')^2}\left(\beta_r^2 r_\phi^2 +\beta_r^2 (r_\phi')^2 + r_\phi^2 (r_\phi')^2\right)}\tilde{\tilde{v}}_n
	\end{aligned}
\end{equation}

The validity of the generic model can be verified by substituting the specific parametrization of a rectangular slider, which is given by
\begin{equation}
	\begin{aligned}
		r_\phi &= \tfrac{b}{2 c_\phi} \\
		c &= \tfrac{b}{2} t_\phi
	\end{aligned}
\end{equation}
It can be verified that the model then collapses onto the model presented in (\ref{eq:dk}) taking into account that for a rectangular slider $\alpha = -\phi$. Furthermore we should take into account the different parametrization that was used for the contact, i.e. $c$ instead of $\phi$.
\begin{equation}
	\dot{c} = \tfrac{b}{2} \tfrac{1}{c_\phi^2} \dot{\phi} = \sqrt{r_\phi^2+(r_\phi')^2}\dot{\phi}
\end{equation}

Clearly the present model is significantly more complex than the model that was retrieved for the rectangular slider. Nonetheless, there are also clear parallels that would suggest the generic model could share the flatness property with the rectangular model. We leave the question of flatness for future work. 
}

\section{Derivation of $\beta_2$}\label{sec:derivation-of-beta2}

To calculate the double integral we transfer to polar coordinates, integrate over the first quadrant and multiply by $4$.
\begin{align}
\beta_2 &= \tfrac{1}{A}\int_{\mathcal{A}} \|\vectorstyle{d}\| \text{d}A  \\
&= \tfrac{1}{ab} \int_{\mathcal{A}} r^2 \text{d}r \text{d}\alpha \\
&= \tfrac{4}{ab} \int_0^{\alpha'} \int_0^{\frac{a}{2}\sec \alpha} r^2 \text{d}r \text{d}\alpha + \tfrac{4}{ab} \int_{\alpha'}^{\frac{\pi}{2}} \int_0^{\frac{b}{2}\csc \alpha} r^2 \text{d}r \text{d}\alpha \\
&= \tfrac{4 a^2}{3 b} \int_0^{\alpha'} \sec^ 3 \alpha \text{d}\alpha + \tfrac{4 b^2}{3 a} \int_{\alpha'}^{\frac{\alpha}{2}} \csc^3 \alpha \text{d}\alpha \\
&= \dots  \\
&=  \frac{a^3 \log \frac{\sqrt{a^2+b^2} + b}{a} - b^3 \log \frac{\sqrt{a^2+b^2} - a}{b} + 2 \sqrt{a^2+b^2} a b }{12 ab}
\end{align}
where $\alpha' = \arctan \frac{b}{a}$.

\section{Reciprocal pressure distributions}\label{sec:reciprocal-pressure-distributions}
Since the pressure distribution is speculative, in conclusion we may try and find reciprocal pressure distributions for either modelling strategies so that the associated geometric factors are equal to the geometric factor from the reciprocal modelling approach. 

To do so we rely on the principle of Maximum Entropy. Here a distribution is sought with maximal entropy but so that it meets additional moment constraints. For details, we refer to \cite{murphy2012machine}.
\begin{equation}
\max_{\rho} \int_{\mathcal{X}}\rho(\vectorstyle{x}) \log \rho(\vectorstyle{x}) \text{d}\vectorstyle{x}
\text{ s.t. } \int f_i(\vectorstyle{x})\rho(\vectorstyle{x}) = \mu_i, ~\forall i
\end{equation}

\begin{enumerate}
\item For the first strategy we have the following problem
\begin{equation}
	\label{eq:ME1}
	\begin{aligned}
		\max_{\rho_1} \int_{\mathcal{A}}\rho_1(\vectorstyle{q}) \log \rho_1(\vectorstyle{q}) \text{d}A \text{ s.t. } & \int_{\mathcal{A}} \rho_1(\vectorstyle{q}) \text{d}A = F\\ & \int_{\mathcal{A}} \|\vectorstyle{q}\|^2 \rho_1(\vectorstyle{q}) \text{d}A = F \beta_2^2 
	\end{aligned}
\end{equation}
The variational problem above can be solved for the so-called MaxEnt distribution. It is well-known that the MaxEnt distribution is of the exponential class. Here $\mu$ and $\lambda$ denote Lagrangian multipliers associated to the two constraints.
\begin{equation}
	\rho_1(\vectorstyle{q}) = \mu \exp(-\lambda \|\vectorstyle{q}^2\|)
\end{equation}
Plugging the MaxEnt distribution in the constraints (\ref{eq:ME1}), yields values for $\mu$ and $\lambda$. For a rectangular slider we have
\begin{equation}
	\beta_2^2 = \frac{\frac{1}{2\sqrt{\pi}}\sum_{\alpha\in{a,b}}\erf\left(\frac{\sqrt{\lambda}\alpha}{2}\right)^{-1}e^{-\frac{\lambda \alpha^2}{4}}\sqrt{\lambda}\alpha -1 }{\sqrt{\lambda}^2}
\end{equation}
which yields an explicit expression for $\lambda$ which can be solved numerically. The value of $\mu$ is more straightforward and simply normalizes the distribution to $F$.
\item Analogous to the method above, for the second modelling strategy we have the following problem
\begin{equation}
	\label{eq:ME2}
	\begin{aligned}
		\max_{\rho_1} \int_{\mathcal{A}}\rho_1(\vectorstyle{q}) \log \rho_1(\vectorstyle{q}) \text{d}A \text{ s.t. } & \int_{\mathcal{A}} \rho_1(\vectorstyle{q}) \text{d}A = F\\ & \int_{\mathcal{A}} \|\vectorstyle{q}\| \rho_1(\vectorstyle{q}) \text{d}A = F \beta_1^2 
	\end{aligned}
\end{equation}
The MaxEnt distribution is given by
\begin{equation}
	\rho_2(\vectorstyle{q}) = \mu \exp(-\lambda\| \vectorstyle{q}\|)
\end{equation}
The computation of $\mu$ and $\lambda$ involve some intractable integrals which can only be evaluated numerically. For conciseness they are not included.
\end{enumerate}

\begin{figure}
\centering
\includegraphics[width=\columnwidth]{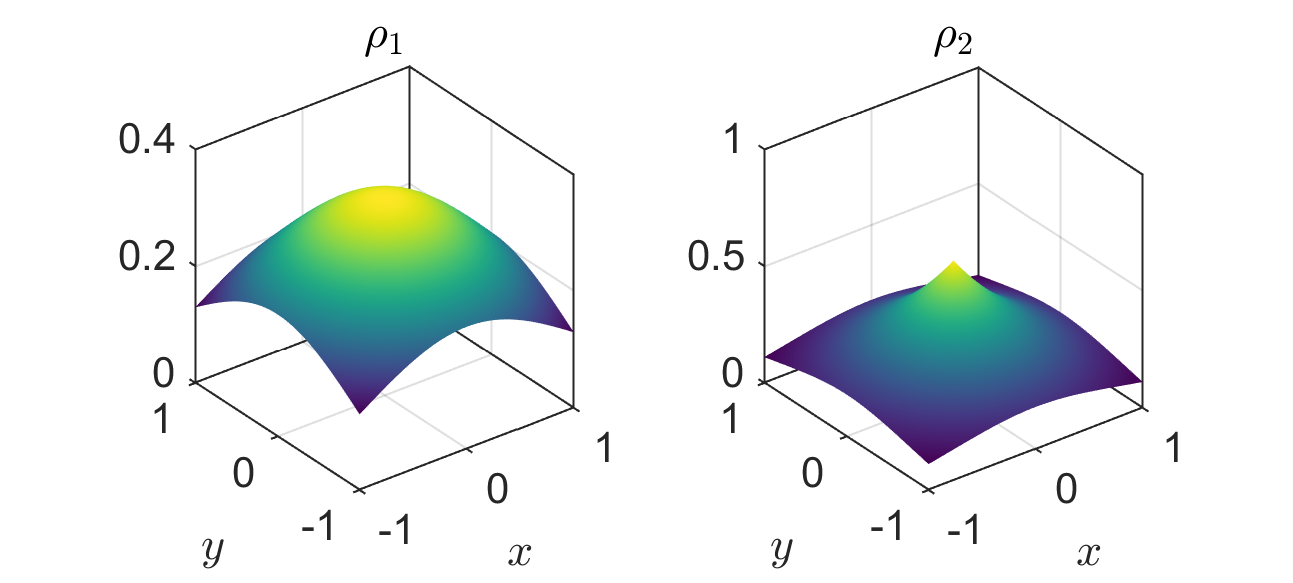}
\caption{Analytical comparison of $\rho_1$ and $\rho_2$.}
\label{fig:rho}
\end{figure}

A comparison between $\rho_1$ and $\rho_2$ is made in Fig. \ref{fig:rho}. Either distribution seems to favour angular motion over linear motion.

\section{Derivation of differential flatness}\label{sec:derivation-of-differential-flatness}
The derivation of $\theta_s$ is straightforward dividing the differential kinematic expressions for $\dot{x}_s$ and $\dot{y}_s$. Differentiating $\theta_s$ yields an expression for $\dot{\theta}_s$
\begin{equation}
\dot{\theta}_s = \tfrac{\dot{x}_s\ddot{y}_s - \ddot{x}_s\dot{y}_s}{\dot{x}_s^2+\dot{y}_s^2}
\end{equation}

By squaring and superposing $\dot{x}_s$ and $\dot{y}_s$ we gain access to the following expression
\begin{equation}
\dot{x}_s^2 + \dot{y}_s^2 = \left(\tfrac{\beta^2}{\beta^2 + c^2}\right)^2 v_n^2
\end{equation}
which we combine with $\dot{\theta}_s$ to derive the flat expression of $c$
\begin{equation}
c = \beta^2 \tfrac{\dot{\theta}_s}{\sqrt{\dot{x}_s^2+\dot{y}_s^2}}
\end{equation}

To find the flat expression for $v_n$ we manipulate the differential kinematic expression for $\dot{\theta}_s$ and substitute the known flat expressions for $C$ and $\dot{\theta}_s$ 
\begin{equation}
v_n = \tfrac{\beta^2 + c^2}{c} \dot{\theta}_s
\end{equation}

Finally we can take the derivative of the flat expression of $c$ and substitute it in its differential kinematic form to arrive at an expression for $v_t$
\begin{align}
\dot{c} &= \beta^2\tfrac{\dot{x}_s\dddot{y}_s - \dddot{x}_s\dot{y}_s}{\sqrt{\dot{x}_s^2+\dot{y}_s^2}^{3}} + 3\beta^2 \tfrac{\left(\ddot{x}_s\dot{y}_s-\dot{x}_s\ddot{y}_s\right)\left(\dot{x}_s\ddot{x}_s+\dot{y}_s\ddot{y}_s\right)}{\sqrt{\dot{x}_s^2+\dot{y}_s^2}^{5}} \\
v_t &= \dot{c} + \left(\tfrac{b}{2}+r\right) \dot{\theta}_s
\end{align}

To connect the differential flat model with that of the kinematic car we proceed as follows. By definition the centre of the kinematic car is given by
\begin{align}
x_p &= x_s + c \cos(\theta_s) + \left(\tfrac{b}{2}+r\right)  \sin(\theta_s) \\
y_p &= y_s + c \sin(\theta_s) - \left(\tfrac{b}{2}+r\right) \cos(\theta_s)
\end{align}
whilst from e.g. \cite{fliess1995flatness} it is known that
\begin{align}
\theta_p &= -\arctan \tfrac{\dot{x}_p}{\dot{y}_p}\\
v &= \sqrt{\dot{x}_p^2+\dot{y}_p^2}\\
\omega &= \tfrac{\dot{x}_p \ddot{y}_p - \ddot{x}_p\dot{y}_p}{\dot{x}_p^2 + \dot{y}_p^2}
\end{align}

Substitution of the expressions for $c$ and $\theta_s$ yields the flat expressions for the car's configuration and input that were stated earlier.

{	\section{Proof of Theorem \ref{th:2}}\label{sec:verification-of-theorem-refth2}

The property is trivially verified for the flat expression for $\theta_s$
\begin{equation}
	\begin{aligned}
		\theta_s &= -\arctan \tfrac{\dot{x}_s}{\dot{y}_s} \\
		&= -\arctan \tfrac{x_s' \dot{\tau}}{{y}'_s\dot{\tau}} \\
		&=  -\arctan \tfrac{x_s'}{y_s'}
	\end{aligned}
\end{equation}

For $c$ the property is less trivial
\begin{equation}
	\begin{aligned}
		c &= \beta^2 \tfrac{\dot{x}_s\ddot{y}_s - \ddot{x}_s\dot{y}_s}{\sqrt{\dot{x}_s^2+\dot{y}_s^2}^{3}} \\
		&= \beta^2 \tfrac{x_s'\dot{\tau}(y_s'\ddot{\tau} + y_s''\dot{\tau}^2) - (x_s'\ddot{\tau} + x_s''\dot{\tau}^2)y_s'\dot{\tau}}{\dot{\tau}^3\sqrt{(x'_s)^2+(y'_s)^2}^{3}} \\
		&= \beta^2 \tfrac{x_s'y_s''-x_s''y_s'}{\sqrt{(x'_s)^2+(y'_s)^2}^{3}} 
	\end{aligned}
\end{equation}
}

\newpage

\bibliography{references}
\end{document}